\documentclass{article}

\usepackage{arxiv}

\usepackage[utf8]{inputenc} 
\usepackage[T1]{fontenc}    
\usepackage{hyperref}       
\usepackage{url}            
\usepackage{booktabs}       
\usepackage{amsfonts}       
\usepackage{nicefrac}       
\usepackage{microtype}      
\usepackage{lipsum}		
\usepackage{graphicx}
\usepackage[round]{natbib}
\usepackage{doi}

\title{A stochastic programming approach for dynamic allocation of bed capacity and assignment of patients to collaborating hospitals during pandemic~outbreaks}

\author{ {Stef Baas}\thanks{Corresponding author} \\
	Department of Stochastic Operations Research\\
	University of Twente\\
	7522 NB Enschede, the Netherlands \\
	\And 
	{Sander Dijkstra} \\
	Department of Stochastic Operations Research\\
University of Twente\\
7522 NB Enschede, the Netherlands \\
	\And
	{Richard J. Boucherie} \\
	Department of Stochastic Operations Research\\
University of Twente\\
7522 NB Enschede, the Netherlands \\
	\And
	{Anne Zander} \\
	Department of Stochastic Operations Research\\
University of Twente\\
7522 NB Enschede, the Netherlands 
}

\renewcommand{\shorttitle}{SP approach for dynamic allocation of bed
	capacity and assignment during
	pandemics}

\hypersetup{
pdftitle={A stochastic programming approach for dynamic allocation of bed capacity and assignment of patients to collaborating hospitals during pandemic~outbreaks},
pdfauthor={Stef Baas, Sander Dijkstra, Richard Boucherie, Anne Zander},
pdfsubject={math.OC},
pdfkeywords={Pandemic response management, Stochastic programming, Simulation, COVID-19, Bed occupancy},
}
\usepackage[utf8]{inputenc}

\usepackage{framed}
\usepackage{mathtools}
\usepackage{tikz}%
\usepackage{multirow}

\usepackage{bbm}

\usepackage{blkarray} %
\usepackage{listings}
\usepackage{verbatim}
\usepackage{a4wide}
\usepackage{diagbox}
\usepackage{import}

\usepackage{setspace}
\usepackage{parskip}
\usepackage{adjustbox}
\usepackage{booktabs} 
\usepackage{float}
\usepackage{tabularx}
\usepackage{xcolor}
\usepackage{subcaption}
\usepackage{graphicx}
\usepackage{placeins}
\usepackage{amsfonts}
\usepackage[utf8]{inputenc}
\usepackage{amsfonts,amssymb,amsthm}
\usepackage{amsmath}
\usepackage{bm}
\usepackage{lpform} 
\usepackage{pdflscape}
\usepackage{rotating}
\usepackage{lineno,hyperref}
\modulolinenumbers[1]

\newcommand*\Bell{\ensuremath{\boldsymbol\ell}}

\usepackage[english]{babel}

\newtheorem{remark}{Remark}

\usepackage{url}

\let\oldcite\cite
\renewcommand{\cite}[2][]{\mbox{\oldcite[#1]{#2}}}
\let\oldcitep\citep
\renewcommand{\citep}[2][]{\mbox{\oldcitep[#1]{#2}}}
\let\oldcitet\citet
\renewcommand{\citet}[2][]{\mbox{\oldcitet[#1]{#2}}}
\begin{document}
\maketitle

\begin{abstract}
Sustaining regular and infectious care during an infectious outbreak requires adequate management support for capacity allocation for regular and infectious patients. During the COVID-19 pandemic, hospitals faced severe challenges, including uncertainty concerning the number of infectious patients needing hospitalization and too little regional cooperation. This led to inefficient usage of healthcare capacity. To better prepare for future pandemics, we have developed a decision support system for central regional decision-making on opening and closing (regular care) hospital rooms for infectious patients and assigning new infectious patients to regional hospitals. Since the relabeling of rooms takes some lead time, we make decisions with a stochastic lookahead approach using stochastic programming with sample average approximation based on scenarios of the number of occupied infectious beds and infectious patients needing hospitalization. The lookahead approach produces high-quality decisions by measuring the impact of current decisions on future costs, such as costs for bed shortages, unused beds for infectious patients, and opening and closing rooms. Our simulation study applied to a COVID-19 scenario in the Netherlands, demonstrates that the stochastic lookahead approach outperforms a deterministic approach as well as other heuristic decision rules such as hospitals acting individually and implementing a pandemic unit, i.e., one hospital is designated to take all regional infectious patients until full. Our approach is very flexible and capable of tuning the model parameters to take into account the characteristics of future, yet unknown, pandemics, and supports sustaining regular care by minimizing the strain of infectious care on the avilable number of beds for regular care.
\end{abstract}

\vbox{{\bf Highlights:}
\begin{itemize}
\item We develop a stochastic direct lookahead approach using two linked stochastic programs to allocate bed capacity to infectious patients and assign infectious patients to hospitals, considering a region of collaborating hospitals during a pandemic. 
\item Our approach differs from current methods in the literature as capacity allocation and patient assignment are decided through a direct lookahead approach on a regional level, taking  into account more factors of uncertainty and the lead time to make complete rooms (not single beds) available and using dynamic estimates based on regional and in-hospital data. 
\item The developed approach is applicable to pandemics or other infectious disease outbreaks. As a demonstration, we consider a case study of the COVID-19 pandemic in a region in the Netherlands.
\item The stochastic direct lookahead approach results in fewer bed shortages for infectious patients while minimizing the average number of hospital rooms utilized by these patients.
\end{itemize}}


\section{Introduction}\label{sect:introduction}

Maintaining consistent care for both regular and infectious patients during an outbreak requires effective management support to allocate capacity appropriately. For example, during the COVID-19 pandemic, we saw a significant rise in demand for healthcare due to infectious patients needing hospital care. At the peak, there were over 150,000 patients hospitalized with COVID-19 in the United States~\citep{ourworldindata2024}. At the same time, the provision of healthcare was constricted by additional safety measures such as (social) distancing and wearing protective gear. At the beginning of the COVID-19 pandemic, there was a lot of uncertainty, and necessary (real-time) data for planning capacity, including platforms to bundle data on a regional or national level, were not yet available. As a consequence, healthcare providers were mainly focused on reserving enough capacity for unforeseen surges of infectious patients and mostly made decisions individually. Regular patients’ surgeries and treatments were postponed, leading to a significant loss of healthy life years~\citep{deGraff2022,Rovers2022}. For example, during COVID-19, cancer diagnoses decreased by 13 to 26 percent, which means that those patients were diagnosed later, in a potentially higher cancer stage, negatively affecting their health outcomes~\citep{Englum2022,Jacob2021}.
On the one hand, when an expected surge in infectious patients did not happen, the reserved capacity was lost, even though it could have been used to serve regular care. On the other hand, in times of overcrowding due to too few reserved capacities or too few capacities in general, infectious patients had to be transferred to other hospitals, sometimes to other regions. Besides being undesirable for the patients who are then far away from loved ones, those transfers required resources such as ambulances and accompanying personnel who could have served other healthcare demand.

We need to prepare for the next pandemic, making use of accurately predicted healthcare demand to then plan healthcare delivery accordingly on a regional level. It is expected that pandemics will again stress our healthcare system. This threat is amplified by the tightening of healthcare capacity due to both an increase in demand, e.g., due to the aging population, and a decrease in available resources, such as a stable decrease of available hospital beds over the last decades~\citep{ Bernstein2022,StatLine2023}. When reserving capacity for infectious care,  note that preparing a regular hospital room for infectious patients requires setup, which takes one to two days (regular patients need to be discharged or transferred to other rooms, the room needs to be cleaned, etc.). Hence, it is important to accurately predict the number of hospital beds used by infectious patients a few days ahead and match it to the corresponding number of reserved rooms for infectious patients. The developments of data platforms during COVID-19 will allow fast access to real-time data in the next pandemic, facilitating accurate predictions~\citep{Baas2021}. However, we still lack a healthcare capacity planning method using those predictions to ensure access to healthcare for both infectious and regular patients during an outbreak on a regional level. 

 In this paper, we propose a joint regional decision-making approach, where we aim to make the best decisions by modeling their effects on the future. This  regional approach is, e.g.,
 in direct support of the managerial decision setting in Dutch regions, where each region has a corresponding consultative body for the acute care chain~\oldcitep{LNAZ}.
 To take advantage of the economies of scale, the regional hospital network must decide as one when and where to open and close rooms for infectious patients and to which hospital to send newly infected patients requiring hospital care. To this end, following the framework and terminology in~\citet{Powell2022} (detailed later in this paper), we develop two linked stochastic \emph{direct lookahead}~(DLA) models by solving two stochastic programs with sample average approximation. The first model decides on opening and closing hospital rooms reserved for infectious care, while the second one decides on the allocation of new infectious patients within the regional network of hospitals. The objective is to minimize the sum of costs for opening and closing rooms, costs for infectious patients that cannot be accommodated, costs for rooms currently being available for infectious care patients, and costs for making rooms ready for infectious care patients.
 The stochastic DLA approach allows us to make the best decisions given the information that we have at that moment, e.g., the current number of hospitalized patients and forecasts of the number of new regional infectious patients who require hospital care and the resulting number of occupied beds in the hospitals for the upcoming days. In addition, the method inherently considers the uncertainties with respect to those numbers, taking into account different possible realizations of the future through scenarios. 

The decision rules deduced from our stochastic DLA approach are compared to two other heuristic decision rules in a simulation study. The two heuristic decision rules mimic the individual decision-making of hospitals, as it might have happened during COVID-19, and the implementation of a pandemic unit, i.e., one hospital is designated to take all regional infectious patients until full. We simulate the evolution of demand for COVID-19 care in hospital wards in the ROAZ (Regional Consultation Body on Acute Care) region~\emph{Acute Zorg Euregio} in the Netherlands at the end of~2021. The numerical results show that our stochastic DLA approach considerably outperforms these heuristic decision rules. Furthermore, we demonstrate that a stochastic lookahead is superior to a deterministic one. Therefore, through joint decision-making supported by our method, regions will be better equipped to face the next outbreak, ensuring regional efficient and fair access to care for all patients.

The research leading to this paper had a direct practical impact, as the comparison of one of these heuristics, the pandemic unit, that was intended to be opened by Acute Zorg Euregio, with regional collaboration led to the conclusion for Acute Zorg Euregio that dynamic regional collaboration would be preferred over a pandemic unit in case of a new pandemic outbreak, see, e.g.,~\citet{artikel_NOS}.

The remainder of this paper is organized as follows. In the next section, a comprehensive literature review is presented. Section~\ref{Sect:model} introduces the sequential decision modeling framework used for capacity allocation and patient assignment and also presents the model to describe the evolution of the infectious bed occupancy in the considered region that is used to generate scenarios. Section~\ref{sect:4-methods} presents the stochastic programs and resulting decision rules for opening rooms and assignment of infectious patients. 
Section~\ref{sect:results} describes the setup and results of a simulation study where the decision rules are compared to two simpler heuristic decision rules with an application to the COVID-19 pandemic. Section~\ref{sect:discussion} concludes the paper and gives ideas for further research.

\section{Literature review}\label{sect:LitRev}

As mentioned in Section~\ref{sect:introduction}, this research proposes central decision-making on bed allocation to infectious patients within a regionally collaborative network of hospitals during a pandemic.   
In the non-pandemic setting, such joint bed capacity allocation to a group of patients within a region was studied in, e.g.,~\citet{litvak2008managing} or \citet{Marquinez2021}. 
During the early weeks of the COVID-19 pandemic, \citet{Kaplan2020} was one of the first to apply a model within an individual hospital to determine~(and allocate) the necessary COVID-19 bed capacity. 
Also, within the setting of an individual hospital,~\citet{Lay2023} developed a discrete-event simulation model for allocating bed capacity to infectious patients. Similar to our approach is that due to the contagiousness of patients, only complete care rooms can be allocated, and these are taken out of regular care processes. 
A similarity between the dynamic programming model for allocating bed capacity to infectious patients in \citet{MA2022} and our approach is that the beds are taken out of the regular care processes, which later in time can be reversed so that beds used for infectious care can be used for regular care again. 
\citet{gokalp2023dynamic} also considered dynamic capacity management of the most crucial resources for COVID-19 care in a single hospital, 
but in their developed \emph{Markov decision process}~(MDP) these resources are not (temporarily) taken away from regular care, but brought into the hospital from an external source. A similarity that \citet{gokalp2023dynamic} shares with our approach is that it takes some lead time to make resources available for the infectious patient group. 

In addition to capacity allocation to infectious patients, our research considers assigning patients to hospitals.
During the COVID-19 pandemic, the assignment of patients to different hospitals was also studied, e.g.,~\citet{Sarkar2021} proposed a data-driven decision-making tool to optimize the assignment of infectious patients,~\citet{Aydin2022} presented a linear programming model and~\citet{Aziz2021} a bi-objective optimization model. 
Unlike our DLA approach, these papers did not make daily patient assignment decisions based on bed occupancy predictions. \citet{Sarkar2021} solved a problem instance where the number of patients to be assigned at that certain point in time follows from a compartmental model, whereas~\citet{Aydin2022} and \citet{Aziz2021} knew for every patient to be assigned the exact time-stamp after which they need hospitalization. 
\citet{dijkstra2023dynamic} developed decision rules for assigning patients to hospitals, aiming at fair balancing of the burden of COVID-19 care on the bed occupancy levels within and across regions. 
The work of~\citet{Ye2022} also presented a patient assignment method based on a defined principle of fairness to balance the occupancy levels in various affected areas. 

Several papers in the literature consider both the allocation of hospital resources and patient assignment during pandemic times.  
Similar to our approach, \citet{parker2020optimal} proposed an optimization model for centralized, daily decisions on patient-hospital assignments and resource re-allocations. 
We conclude these decisions by solving two stochastic programs, while \citet{parker2020optimal} use a fixed-horizon robust optimization program, where the (uncertainty set of) daily number of new patient arrivals is assumed to be within in an interval, and the \emph{length of stay}~(LoS) of each patient follows from a Weibull distribution. Their objective is to minimize the total bed shortage in the network of collaborating hospitals.
\citet{Eriskin2022} proposed a robust optimization model that determines the periodic resource re-allocations and daily patient-hospital assignments. This model is a $p$-robust program that includes stochasticity regarding bed occupancy. 
\citet{Barbato2023} proposed a mathematical programming approach for 
resource re-allocations and re-purposing of hospital wards~(i.e., converting ordinary beds into isolation beds) and, as a last resort, the selective discharge of less severely ill patients. 
In comparison to our approach, the model in \citet{Barbato2023} is not evaluated in a dynamic way, and parameters, such as the number of patients to be assigned, were assumed to be deterministically known. 
\citet{fattahi2023resource} developed a multi-stage stochastic program to make three types of decisions: redistribution of patients between hospitals; resource allocation, i.e., providing more external resources for hospitals; and resource relocation. Their objective was to minimize patient refusals, the added external resources, and the number of both patient and resource transfers. Their resulting policy prescribed decisions for a fixed-length horizon and was evaluated in a rolling horizon fashion via a simulation model. 

Prior to the COVID-19 pandemic outbreak,  
\citet{Sun2014} was, to the best of our knowledge, the only research that developed multi-objective optimization models for both patient and resource allocation among hospitals during an influenza outbreak. 
Opposed to our approach, where we take capacity allocation and patient assignment decisions every day, they develop an optimization model for a fixed planning horizon.  
Furthermore, their objectives focus on total and maximum patient traveling distance instead of capacity utilization, and their approach does not incorporate stochasticity, i.e., the LoS is assumed to be known for every patient, which also holds for the number of patients to be assigned (before solving the optimization model, the number of new patients to be assigned is known for each day of the seven-days planning horizon). 

In our approach, we make daily decisions on centrally allocating bed capacity to infectious patients in all of a region's hospitals and patient-to-hospital assignments. 
The structure we deploy, in which the optimal solution of a first stochastic program for capacity allocation is input for a second stochastic program that decides upon patient assignment, is new. 
Moreover, both capacity allocation and patient assignment are done dynamically (on a daily basis), and therefore, these decisions are reconsidered or altered after only a short amount of time, contrary to the literature, where assignment decisions are usually made at once for a fixed time period. 
The first part of our model, i.e., the stochastic program to decide upon capacity allocation, has several important components: only complete care rooms can be allocated to infectious patients, it takes a lead time to make this capacity available, and the allocated beds are taken out of regular care processes, which later in time can be reversed. As we described, we have seen models in the literature in the setting of a single hospital that incorporate some of these components. Still, we have not seen a model that combines all three of them (either in a regional setting or in a single-hospital setting). In that sense, we contribute to the literature by modeling a novel problem situation that incorporates all these three components in the capacity allocation part. 
Our stochastic DLA approach allows us to make the best decisions given the information that we have at that moment, taking into account different possible scenarios of the future capacity demand in both stochastic programs. Our approach to generating these scenarios adds another modeling method of the involved stochasticity to the literature, namely one where scenarios of bed occupancy result from inflow prediction and estimation of the LoS for each hospital. 
Compared to studies that take stochasticity into account in their optimization models, the degree to which we include the uncertainty in the underlying bed occupancy in each hospital was not observed in previous work. Our approach takes more factors of uncertainty into account, and the way scenarios are generated is more advanced and accurate. 
A final, short remark on our results is that we provide an extensive simulation study to show the performance of our approach in an environment that includes randomness. This is different from what the majority of literature shows, namely the evaluation of a method in an instance-based manner.

\section{Model, optimization problem, and proposed solution approach} \label{Sect:model}
In this section, we model collaborative bed allocation and patient assignment as a sequential decision problem and motivate our stochastic \emph{direct lookahead}~(DLA) solution approach to that problem. Section \ref{modelframework} explains the general setting of regionally collaborating hospitals~(Section~\ref{model:region}), how we model bed occupancy within a single hospital~(Section~\ref{sect:ModOcc}), and the
forecasting approach for infectious patient arrivals to the hospital and hospital bed occupancy, which will then be used to generate scenarios for the stochastic programs (Section~\ref{sect:ScenGen}). Section~\ref{prob_soln} presents the optimization problem~\ref{sec: SDMP} as well as the solution approach~\ref{sect:PowellPolicies}, further detailed in Section \ref{sect:4-methods} through the exact formulation of the stochastic programs.

\subsection{Modeling framework} \label{modelframework}
\subsubsection{Region of collaborating hospitals} \label{model:region}
We consider a region of $H$~collaborating hospitals, indexed~$h=1,\ldots, H$, that centrally and dynamically decide on their infectious bed capacity and assignment of infectious patients to those hospitals during an infectious outbreak. This collaboration aims to ensure enough capacity for infectious patients while maintaining regular care to the highest possible degree. 

Each hospital has a dedicated ward for infectious patients with bed capacity~$c_h$. In addition, there are $n_h$~regular care rooms at hospital~$h$ that may be opened for infectious care. Hence, infectious bed capacity can only be made available in terms of rooms and not in terms of individual beds for health safety reasons. These rooms in hospital~$h$ are labeled~$1,\ldots,n_h$ and can only be opened and closed in sequence: room $n$ can be opened only if room $n-1$ is open, and room $n$ can be closed only if room $n+1$ is closed. The reason for this is that there must be a clear separation between infectious care rooms and regular care rooms to lower the risk of disease spread. Let room~$n$ in hospital~$h$ contain $b_{h,n}$~beds. In accordance with practice, we assume that it takes a fixed amount of two~days to empty a regular care room and open it for infectious patients, e.g., by relocating regular patients to other regular care rooms and/or discharging regular patients. If a room allocated to infectious patients is empty, then it may be closed and will be immediately available for regular patients. Further, over the course of a day, infectious patients arrive either autonomously at individual hospitals or in the region, in which case they have to be assigned to individual hospitals.


\subsubsection{Modeling bed occupancy}\label{sect:ModOcc}
We assume that the number of infectious patients demanding care follows an inhomogeneous Poisson process with rate~$\lambda_{\tau}$ at time $\tau$. The assumption of Poisson arrivals was shown to be justified for arrivals to Emergency Departments~\citep{WhittSEHpaper} and is also verified for arrivals of COVID-19 patients \citep{dijkstra2023dynamic}. A fraction~$f_h$ of the regional patients arrives autonomously at hospital~$h$, leading to an autonomous inflow with Poisson rate~$ \lambda_{h,\tau}=f_h\lambda_{\tau}$ of infectious patients to hospital~$h$. This autonomous inflow contains, e.g.,  patients who are diagnosed positive upon arrival at the hospital's emergency department. The remaining inflow of regional patients follows a Poisson process with rate~$(\lambda_{\tau}-\sum_h \lambda_{h,\tau})$. The assignment of those patients to hospitals will be part of our solution approach.

Following~\citet{Baas2021}, we model the infectious ward of each hospital $h$ as an infinite server queue that records the number of hospitalized infectious patients in the ward.  Hence, each hospital is described by an~$M_{\tau}/G/\infty$ queue. Let $L_h$ denote the random variable of the \emph{length of stay}~(LoS) of the patients in the ward. Following~\oldcitet[Theorem 1.2]{MasseyWhitt}, starting from an empty system at time 0, the occupancy by autonomously arriving infectious patients~$N'_{h,{\tau}}$ in hospital~$h$ at time~${\tau}$  has a time-dependent Poisson distribution, i.e.,
\begin{equation}
\label{eq:Poissontotal}
\mathbb{P}[N'_{h,{\tau}}=n] = \frac{(\rho_{h,{\tau}})^{n}}{n!} {\rm e}^{-\rho_{h,{\tau}}},
\mbox{ where }
\rho_{h,{\tau}} = \mathbb{E}\left[\,\int^{\tau}_{u=\max\{0,{\tau}-L_h\} }\lambda_{h,u} du\right].
\end{equation}
The  Poisson distribution for the number of (autonomously arriving) infectious patients in the ward~\eqref{eq:Poissontotal} allows us to evaluate various performance measures.
Let~${\bf L}_{h,{\tau}}$ denote tuples of the  attained LoSs (up to time ${\tau}$) of patients residing in the ward of hospital~$h$ at time~${\tau}$. The expected occupancy (only considering the autonomous inflow) in the ward at time~${\tau}+\sigma$ given the LoSs of the residing patients at time~${\tau}$ is determined by the patients present at time~${\tau}$ who are still present at time~${\tau}+\sigma$ and the patients autonomously arriving between time~${\tau}$ and time~${\tau}+\sigma$ in the $M_{\tau}/G/\infty$ queue that starts empty at time~${\tau}$:
\begin{equation}
\label{eq:Poissontotalcond}
\mathbb{E}[N'_{h,{\tau}+\sigma}\;|\;{\bf L}_{h,{\tau}}= \Bell] =  \sum_{i=1}^{| \Bell| }\frac{1-F_h(\ell_i+\sigma)}{1-F_h(\ell_i)}  + \mathbb{E}\left[\,\int^{\tau+\sigma}_{u=\max\{\tau,{\tau+\sigma}-L_h\}}\lambda_{h,u}\, du\right].
\end{equation}
Note that, in~\eqref{eq:Poissontotalcond}, we may include regional patients that have been assigned to hospital~$h$ before time ~${\tau}$ on top of the autonomously arriving patients. Other measures, such as the variance and quantiles of the occupancy at time~${\tau}+\sigma$ may be calculated from similar principles~\oldcitep[see, e.g.,][]{Baas2021}. 

\subsubsection{Scenario generation}\label{sect:ScenGen}

The stochastic programs that we present in Section~\ref{sect:4-methods} use scenarios of daily occupancy and infectious patients arriving in the region to be assigned to the hospitals. These scenarios are generated based on an estimate of the parameters of the model described in Section~\ref{sect:ModOcc}. 
Following the approach developed in \citet{Baas2021} and \citet{dijkstra2023dynamic}, using historical data up to day $t$, we use a prediction $\hat{\lambda}_{\tau}$ for $\tau\in[t,t+s)$ of the regional arrival rate in the period~$[t,t+s)$, an estimate $\hat{f}_h$ of the fraction of regional patients arriving at each hospital $h$, and the (known) LoS distribution. 
Here,~$s$ is the prediction horizon in days, which is usually not more than a week. 
Combining $\hat{\lambda}_{h,\tau}$ for $\tau\in[t,t+s)$ with the LoS distribution, a sample path $\tilde{N}_{h,t+1}, \tilde{N}_{h,t+2}, \dots, \tilde{N}_{h,t+s}$  of the occupancy at days~$t+1,\dots,t+s$ by currently residing and autonomously arriving infectious patients at hospital~$h$ is generated. This method was seen to yield accurate forecasts of daily bed occupancy for the COVID-19 pandemic \citep{Baas2021, dijkstra2023dynamic}. 

Following Section~\ref{sect:ModOcc}, letting $\partial\hat{\Lambda}_{t+u}=\int_{t+u-1}^{t+u}\hat{\lambda}_\tau d\tau$, the autonomous arrival rate of regional infectious patients to hospital~$h$ in~$[t+u-1,t+u)$ can be predicted as $\partial\hat{\Lambda}_{h,t+u} =\hat{f}_h\partial\hat{\Lambda}_{t+u}$ for all $u\in\{1,\dots, s\}$. 
Using this predictor, scenarios $A_{t+u}$ of daily regionally arriving infectious patients to be assigned to hospitals in $[t+u-1,t+u)$ can be determined for~$u=1,\dots,s$ by generating sample paths from the inhomogeneous Poisson process with \mbox{intensity~$((1-\sum_h\hat{f}_h)\partial\hat{\Lambda}_{t+u})_{u=1}^{s}.$ }

\begin{remark}[Poisson arrivals]
In this paper, it is assumed that the arrival process of regional infectious patients is a nonhomogeneous Poisson process. 
 This assumption was verified for arrivals to Emergency Departments in~\citet{WhittSEHpaper}, while it resulted in an accurate forecast of occupancy by COVID-19 patients in~\citet{Baas2021} and~\oldcitet{dijkstra2023dynamic}. The assumption of Poisson arrivals is not strictly necessary for our proposed approach, which only uses scenarios of the number of patients that arrived on a given day. 
  For the proposed approach, it is necessary that one can reliably estimate and predict  the distribution of the daily arrivals based on historical arrival data. 
 When the assumption that the arrivals follow an inhomogeneous Poisson process is relaxed, the formulas in Section~\ref{Sect:model} are not guaranteed to hold anymore, however, and hence it may be more difficult to explain certain behavior of operating characteristics of the model.  
\end{remark}

\subsection{Problem and solution approach}\label{prob_soln}
\subsubsection{Sequential decision-making problem} \label{sec: SDMP}
We model the capacity allocation and patient assignment as a sequential decision problem according to the unified framework proposed in \citet{Powell2022}. The \emph{state} variables indicate the rooms that are currently open (in the morning) on day~$t$ and that are in preparation for opening on day~$t+1$ and~$t+2$ in hospital~$h$, as well as the current number of infectious patients per hospital and their attained LoS. The \emph{decisions} are which rooms for each hospital~$h$ will actually be opened on day~$t$ (which had to be prepared on day~$t-2$ and~$t-1$ to be opened on day~$t$), which rooms will be prepared to be opened on day~$t+1$ (which had to be prepared on day~$t-1$ and~$t$ to be opened on day~$t+1$), which rooms will be prepared to be opened on day~$t+2$ and which rooms to close. Further, the hospitals decide online which regional patient to assign to which hospital. In between decision epochs, the \emph{exogenous information} becomes known, i.e., the number of autonomously arriving infectious patients and the number of discharged patients. The \emph{transition} to the next state happens through implementing the decisions and taking the exogenous information into account. Our \emph{objective} is to minimize the cumulative daily costs, which consist of costs for opening and closing rooms (reflecting cleaning and set-up costs), costs for overbeds (infectious patients that cannot be accommodated), costs for the number of reserved infectious beds (beds in rooms open for infectious patients) and costs for the number of infectious beds that are in the opening process (beds in rooms from which regular patients are transferred or discharged).

\vbox{\begin{remark}[Overbeds]
    In this paper, it is assumed that when the infectious occupancy at a hospital exceeds the capacity, the remaining patients are placed on an extra bed, termed an overbed, which incurs an additional cost to the hospital and which may be realized by temporarily reducing the nurse-to-patient ratio. 
    The term overbed was also used in~\citet{Baas2021, dijkstra2023dynamic}.
    The assumption that these patients can always be placed on overbeds is only valid when the expected number of overbeds per day is small. In situations where there is not enough capacity in the region or when there is a very sudden and high increase in infectious patients that need admission, the quality of the solution provided by our approach is no longer guaranteed. In such cases, the proposed approach could be combined with other approaches that deal with assigning patients to other regions, such as the one proposed in~\citet{dijkstra2023dynamic}. 
\end{remark}}

\subsubsection{Solution approach}\label{sect:PowellPolicies}
The two time-scale sequential decision problem of opening and closing rooms on a daily basis and assigning patients to hospitals during the day corresponds to an infinite horizon MDP with large discrete state and action spaces and partially unknown dynamics. 
In particular, it is assumed that the regional arrival rate $\lambda_\tau$ and the fractions $f_h$ of autonomous regionally arriving infectious patients to each hospital are unknown, while the distribution of $L_h$ is known. However, we do assume that we can reasonably predict the regional arrival rate over a short time horizon of a few days based on the arrival history. 
Therefore, the MDP cannot directly be solved to optimality for real-life instances. Instead, we aim to find a good approximate policy, i.e., a mapping from states to actions. An approximate policy can be constructed based on one or a mix of the four meta-policies for solving sequential decision-making problems~\citep{Powell2022}. Since our focus is on making good decisions for the current decision epoch, we use a (stochastic) DLA approach, which estimates the impact of current decisions on current and future costs. Through modeling relationships between succeeding decision epochs, we can avoid the introduction of parameters that are required for the other meta-policies. This means that we can directly apply our DLA method to make decisions without prior computation to tune parameters. 

In reality, a (stochastic) DLA model should be solved to include all up-to-date information at every decision time. That applies, in particular, to decisions for arriving regional infectious patients who need to be assigned to a hospital throughout the day. To ensure manageable run times during the evaluation phase of our approach via simulation, we decide to take all patient steering decisions for one day at once at the beginning of a day by deciding to which hospital to assign the $n$th patient of the day. Once actual patients arrive throughout the day, they are assigned to hospitals based on those decisions. However, this model can easily be adapted to make online decisions on an individual patient level (to be applied in the real world).
Therefore, we develop two linked stochastic DLA models for the two types of decisions~(on opening and closing of rooms and assigning regional infectious patients) that we need to make. To build tractable lookahead models, we truncate the time horizon (looking at least two days ahead), we sample scenarios with respect to the exogenous information, and we apply stage aggregation, in this case working with two stages. Because we are working with discrete and constraint decisions and assume linear costs, we develop two 2-stage stochastic programs, which we apply in a rolling horizon fashion. 
The stochastic programs are explained in more detail in Section~\ref{sect:4-methods}.

\section{Decision rules for opening rooms and patient assignment} \label{sect:4-methods}

In this section, we will introduce the two stochastic programs and the resulting decision rules. The program~SP1, which decides on the number of regular care rooms that are made available and planned to be made available for infectious patients in the coming days, is described in Section~\ref{Subsect: ILP_opening_closing}, while the program SP2, which decides on how to assign patients which will arrive throughout the current day to hospitals, conditional on the previous decision to open rooms for that day made by SP1, is given in Section~\ref{Subsect:patient_allocation}. Section~\ref{sect:relationship_SPs} comments on the relationship between~SP1 and~SP2.

\subsection{Stochastic program for room allocation} \label{Subsect: ILP_opening_closing}

This section introduces the stochastic program SP1 for opening and closing rooms for infectious care patients.  In the first stage, SP1 decides on the number of regular care rooms that are actually opened on day~$t$ and prepared to be made available on days $t+1, t+2$  for infectious patients. In the second stage of SP1, we assign future regional infectious patients to the hospitals, which are, hence, wait-and-see variables. The objective is to minimize a weighted sum of the number of overbeds, available and reserved regular care beds, and opened and closed regular care rooms.  

\subsubsection{Objective}
SP1 minimizes the average daily costs over the days until the end of the considered planning horizon, where the average is taken over scenarios.
The objective function equals a linear combination of the number of opened and closed rooms, regular care beds used and scheduled to be used, as well as the number of overbeds, with coefficients $\alpha, \beta, \gamma, \delta, \epsilon>0$, respectively.
Note that all costs besides those for overbeds can be incurred at the time of decision-making. The costs for overbeds can only be incurred the next day when the occupancies become known. 
At the end of the considered lookahead horizon, the remaining number of overbeds and opened rooms are multiplied by costs $(s-1)\cdot \epsilon$~and~$(s-1)\cdot \gamma$, respectively, and added to the objective function as a terminal cost so that the final amount of overbeds and opened rooms are counted for $s$ days in total. 
The choice of the terminal cost was made to incentivize consistency in the decisions made by the stochastic program. 
Suppose the terminal cost term is not added. In that case, the program can choose to allow for a large number of overbeds on day~$t+s$ or opened rooms on day~$t+s-1$ when running the program at day $t$, but when the program is actually used to decide on open rooms on day~$t+s-1$ these overbeds or used beds induce a higher cost and will be avoided as a result.

\subsubsection{Decision variables and scenarios} The primary decision variables for SP1 are indicator variables~$(z_{h,n,t+u})_{u=0}^{s-1}$, each denoting whether room~$n$ in hospital~$h$ is opened at a day~$t+u$. As secondary decision variables, we determine indicator variables $(v_{h,n,t+u}^{(d)})_{u=0}^{s-1}$, each denoting whether room $n$ in hospital $h$ is prepared to be open on day~$t+u+d$ at a day $t+u$. 
As it was assumed in Section~\ref{Sect:model} that a regular care room can only be opened for infectious care after 2~days, given that it is currently in use for regular care, we must have that $z_{h,n,t+u}\leq v_{h,n,t+u-1}^{(1)}\leq v_{h,n,t+u-2}^{(2)}$, i.e., rooms opened at a current day were scheduled to be opened yesterday and the day before yesterday. 
A third set of decision variables are non-negative integers~$(x^{(i)}_{h,t+u})_{u=1}^{s}$,  each denoting the number of regionally arriving infectious patients assigned to hospital $h$ in a period $[t+u-1,t+u)$, depending on scenario $i$,
 which is defined as an independently sampled path~$((\tilde{N}^{(i)}_{h,t+u})_{h=1}^H,A_{t+u}^{(i)})_{u=1}^s$ 
of the occupancy at the hospitals and regionally arriving infectious patients.
Based on these three sets of decision variables and the scenarios, we can determine the offered capacity at days $t,\dots, t+s-1$, and corresponding occupancy at (the beginning) of the days $t+1,\dots,t+s$. 
In this calculation, it is assumed that regional patients assigned to a hospital stay at the hospital for each day in the time interval $[t,t+s)$. The difference between occupancy and capacity, if it is positive, results in the number of overbeds for that day~\oldcitep[see also][]{dijkstra2023dynamic}. 

\subsubsection{Implementation and reuse of first-stage decisions}
The first-stage decision variables~$(z_{h,n,t+u})^{s-1}_{u=0},\, v_{h,n,t}^{(1)}$, and~$v_{h,n,t}^{(2)}$ are stored after determining the optimal solution. We implement the decisions for $u=0$ and take the decision variables~$(z_{h,n,t+u})^{s-1}_{u=0}$ into account in SP2, i.e., when assigning patients to hospitals. Further, the decision variables~$z_{h,n,t}, v_{h,n,t}^{(1)}$ and $v_{h,n,t}^{(2)}$ are used in the first stochastic program when running it to determine rooms to open on the next day. 

 \subsubsection{Stochastic program}
Table~\ref{table: model definitions2} describes the indices, parameters and optimization variables used in~SP1, where we define $\mathcal{U}_-=\{0,\dots, s-1\}$ and $\mathcal{U}_+=\{1,\dots, s\}$.

The stochastic program SP1 for opening and closing rooms is defined below. Unless indicated otherwise, the index variables are assumed to lie in the domains indicated in Table~\ref{table: model definitions2}.

\renewcommand{\lpforall}[1]{&&\hspace{-40mm} \forall #1}
\begin{lpformulation}[\text{(SP1)}]
\lplabel{ILPX}
\lpobj[ILP_shorter obj]{min}{\textstyle\sum_{u=0}^{s-1}\sum_{h=1}^H[\alpha y_{h,t+u}^{(+)} +\beta y_{h,t+u}^{(-)}+  \sum_{n=1}^{n_h} (\gamma z_{h,n,t+u}
+\delta ( v_{h,n,t+u}^{(1)} +  v_{h,n,t+u}^{(2)} ))  b_{h,n}]\lpnewline +\frac{1}{I}\textstyle\sum_{i=1}^I\sum_{h=1}^H[ (\epsilon\sum_{u=1}^so_{h,t+u}^{(i)})+ (s-1)(\epsilon o_{h,t+s}^{(i)} +\gamma z_{h,n,t+s-1} )]}
\lpeq*{\text{\it First stage:}}{}{}
\lpeq[ordering open_close]{z_{h,n,t+u}\leq z_{h,n-1,t+u}}{u\in\mathcal{U}_-,n\geq 2,h,}
\lpeq[determine opened_closed rooms]{\textstyle\sum_{n=1}^{n_h}(z_{h,n,t+u} - z_{h,n,t+u-1}) = y_{h,t+u}^{(+)} - y_{h,t+u}^{(-)}}{u\in\mathcal{U}_-,h,}
\lpeq[only_open_if_reserved]{z_{h,n,t+u}\leq v^{(1)}_{h,n,t+u-1}+ z_{h,n,t+u-1}}{u\in\mathcal{U}_-,h,n,}
\lpeq[only_res_tomorrow_if_dayafter]{v^{(1)}_{h,n,t+u}\leq v^{(2)}_{h,n,t+u-1} +v^{(1)}_{h,n,t+u-1}  }{u\in\mathcal{U}_-,h,n,}
\lpeq[closed_yest_not_open_tom]{z_{h,n,t+u+1} - z_{h,n,t+u}-1\leq z_{h,n,t+u} - z_{h,n,t+u-1}}{u\in\mathcal{U}_-,h,n,}
\lpeq[binary vars]{ z_{h,n,t+u}, \,v_{h,n,t+u}^{(1)},\,v_{h,n,t+u}^{(2)}\in\{0,1\}}{u\in\mathcal{U}_-,h,n,}\\\nonumber
\lpeq*{\text{\it Second stage:}}{}{}
\lpeq[regional allocation]{\textstyle\sum_{h=1}^H x_{h,t+u}^{(i)}=A_{t+u}^{(i)}}{u\in\mathcal{U}_+, i,}
\lpeq[determine occupancy]{\textstyle N_{h,t+u}^{(i)}=\tilde{N}_{h,t+u}^{(i)}+ \sum_{\ell=1}^{u} x_{h,t+\ell}^{(i)}}{u\in\mathcal{U}_+,h,i,}
\lpeq[determine overbeds]{N_{h,t+u}^{(i)}-\textstyle\sum_{n=1}^{n_h}z_{h,n,t+u-1}b_{h,n} -c_h\leq o_{h,t+u}^{(i)} }{u\in\mathcal{U}_+,h,i,}
\lpeq[natural vars]{x_{h,t+u}^{(i)}, \,N^{(i)}_{h,t+u},\,o_{h,t+u}^{(i)} \in\mathbb{N}_0}{u\in\mathcal{U}_+,h,i,}
\lpeq[natural vars2]{y_{h,t+u}^{(+)},\,y_{h,t+u}^{(-)} \in\mathbb{N}_0}{u\in\mathcal{U}_-,h.}
\end{lpformulation}
\FloatBarrier

In the above problem formulation, the Constraint set~\eqref{ordering open_close} enforces that rooms can only be opened sequentially.
Constraint set~\eqref{determine opened_closed rooms} determines the number of regular care rooms opened or closed at a hospital on a given day.
 Constraint set~\eqref{only_open_if_reserved} ensures that rooms can only be open on a given day if that room was open the day before or if the room was scheduled to be open on that day the day before. 
Constraint  set~\eqref{only_res_tomorrow_if_dayafter} ensures that a room can only be scheduled to open tomorrow if it was already scheduled to open in two days the day before~(which ensures that opening a room takes two days) or if the room was scheduled to be open in one day the day before.
 Constraint set~\eqref{closed_yest_not_open_tom} ensures that rooms that are closed on a given day cannot be opened the next day. Constraint set~\eqref{binary vars} considers domain constraints.
Constraint set~\eqref{regional allocation} ensures that all regionally arriving infectious patients in need of assignment are assigned to a hospital for all scenarios and considered days.  Constraint set~\eqref{determine occupancy} determines the total occupancy at a hospital based on the assigned regional patients and the autonomous occupancy under the scenario on a given day. 
Constraint set~\eqref{determine overbeds} determines the amount of overbeds at a hospital under a scenario on a given day.
 Constraint sets~\eqref{natural vars} and~\eqref{natural vars2} are domain constraints.

\begin{table}[h]
\small
\caption{Indices, parameters, and decision variables used in the stochastic program SP1 for opening and closing regular care rooms. The function $\mathbb{I}(\cdot)$ denotes the indicator function.}
\renewcommand{\arraystretch}{1.2}
\centering
\begin{tabular}{@{} l l @{}}
\toprule
\textbf{Symbol} & \textbf{Description} \\
[0.5ex]
\midrule
\addlinespace[1ex]
\textbf{Indices}  & \\
$h\in\{1,\dots, H\}$ & Hospital\\
$n\in\{1,\dots, n_h\}$ & Room in hospital $h$\\
$u\in\{0,\dots, s\}$ & Days ahead  \\
$i\in\{1,\dots,I\}$ & Scenario (scen.) \\
$d\in\{1,2\}$ & Day ahead in schedule\\
\addlinespace[1ex]
\textbf{Parameters} & \\
{\it First stage:}\\
$b_{h,n}\in\mathbb{N}$ & Number of beds in room $n$, hospital $h$\\
$z_{h,n,t-1}\in\{0,1\}$ & $\mathbb{I}($room $n$ in hospital $h$ is open on day~$t-1)$ \\
$v_{h,n,t-1}^{(d)}\in\{0,1\}$& $\mathbb{I}($room $n$, hospital $h$ scheduled open in $d$ days on day~$t-1)$\\
$c_h\in\mathbb{N}$& Standard capacity infectious ward in hospital $h$\\
$\alpha, \beta,\gamma, \delta,\epsilon\in\mathbb{R}_{>0}$ & Weights used in objective function\\
$M\in\mathbb{R}_{>0}$  & Big number, set to~$10^7$\\\\
{\it Second stage:}\\
$\tilde{N}^{(i)}_{h,t+u}\in\mathbb{N}$ & Autonomous  occupancy, hospital $h$, on day~$t+u$, $u\in\mathcal{U}_+,$  scen. $i$ \\
    $A_{t+u}^{(i)}\in\mathbb{N}$   & Regional patients arriving
    in $[t+u-1,t+u)$, $u\in\mathcal{U}_+$, scen. $i$\\
\addlinespace[1ex]
\textbf{Variables} & \\
{\it First stage:}\\
$z_{h,n,t+u}\in\{0,1\}$ & $\mathbb{I}($room $n$ is open  at hospital $h $ on day~$t+u$), $u \in\mathcal{U}_-$\\
$v_{h,n,t+u}^{(d)}\in\{0,1\}$ & $\mathbb{I}($room $n$, hospital $h$ scheduled open in $d$ days on day~$t+u)$, $u \in\mathcal{U}_-$\\
$y_{h,t+u}^{(+)}\in\mathbb{N}_0$ &  Total infectious rooms opened  in hospital $h$ on day~$t+u$, $u\in\mathcal{U}_-$\\
$y_{h,t+u}^{(-)}\in\mathbb{N}_0$ & Total infectious rooms closed in hospital $h$ on day~$t+u$, $u\in\mathcal{U}_-$\\
\\
{\it Second stage:}\\
$x_{h,t+u}^{(i)}\in\mathbb{N}_0$ & Patients assigned to hospital $h$ in $[t+u-1,t+u)$, $u\in\mathcal{U}_+$, scen. $i$\\
$N_{h,t+u}^{(i)}\in\mathbb{N}_0$ & Total occupancy in hospital $h$ on day~$t+u$, $u\in\mathcal{U}_+,$ scen. $i$\\
$o_{h,t+u}^{(i)}\in\mathbb{N}_0$ & Overbeds in hospital $h$ on day~$t+u$, $u\in\mathcal{U}_+,$ scen. $i$\\
\bottomrule
\end{tabular}
\label{table: model definitions2}
\end{table}
\FloatBarrier



\subsection{Stochastic program for patient assignment~\label{Subsect:patient_allocation}}
This section introduces the stochastic program SP2 for the assignment of patients to hospitals during the current day. It takes the first stage decisions on rooms of SP1 as input. The first stage of SP2 decides on the assignment of the~$j$-th regional patient of today for $j=1,2,\dots$ simultaneously, and the second stage decides on the assignment of future regional patients. Hence, we solve SP2 once directly after SP1 to decide on the assignment of all future regional infectious patients of the day. The objective of the program is to minimize the number of overbeds for infectious patients in the coming days $t+1,\dots, t+s$. 

\subsubsection{Objective}
The stochastic program minimizes the objective function, which equals the average number of overbeds over the days~$t+1,\dots,t+s$ that are related to decisions made on the \hbox{days~$t,\dots, t+s-1$,} where the average is taken over the chosen scenarios.

\subsubsection{Decision variables and scenarios}
The probability $p^{(j)}$ that $j$ regional patients need to be assigned today coming from the Poisson distribution with rate~$\partial\tilde{\Lambda}_{t+1} = (1-\sum_h\hat{f}_h)\partial\hat{\Lambda}_{t+1}$ is given as an input parameter to~SP2. 
These probabilities are truncated at a value~$J$, set to the $97.5\%$~quantile of the Poisson distribution with rate~$\partial\tilde{\Lambda}_{t+1}$.
Similar to the scenarios for~SP1, a scenario  $(i,j)$ for~SP2 is defined as an independent sample path~$(((\tilde{N}^{(i)}_{h,t+u})_{h=1}^H)_{u=1}^s, (A_{t+u}^{(i)})_{u=2}^s,j)$ of the autonomous occupancy at the hospitals at days $t+1,\dots, t+s$, regionally arriving infectious patients in need of assignment in $[t+1,t+s)$ and where $j$ is the number of infectious patients in need of assignment in $[t,t+1)$.

The primary decision variables of the stochastic program for patient assignment come in the form of an indicator variable $w^{(j)}_{h,t+1}$ of the event that the~$j$-th arriving patient in the period $[t,t+1)$ will be assigned to hospital~$h$.
A second set of decision variables $(x^{(i,j)}_{h,t+u})_{u=2}^s$  represents the number of regionally arriving infectious patients assigned to hospital $h$ in $[t+u-1,t+u)$ for scenario~$(i,j)$.
Given the patient assignments and the autonomous occupancy scenario, the total occupancy can be determined.  The capacity $C_{h,t+u}$ at the hospitals at time $t+u$, determined as \hbox{$C_{h,t+u}=c_h+\sum_{n=1}^{n_h}z_{h,n,t+u}b_{h,n}$,} is determined by the solution of~SP1. The capacity over time and occupancy scenarios now determine the number of overbeds $(o^{(i,j)}_{h,t+u})_{u=1}^s$ for each scenario~$(i,j)$, hospital $h$, and time $t+u$. 

\subsubsection{Implementation of first-stage decisions}
The first-stage stage decision variables~$w^{(j)}_{h,t}$ are stored after solving SP2 to optimality. When new regional infectious patients arrive throughout the day, we assign them to hospitals according to the values of those first-stage variables. If the actual number of patients that arrive during a given day exceeds the upper bound~$J$, the choice is made to run the stochastic program again with limit~$2\cdot J$ setting~$p^{(j)}=0$ for~$j\leq J$, in order to obtain decision variables for a larger amount of arriving patients.

\subsubsection{Stochastic program}
Table~\ref{table: model definitions3} describes the indices, parameters, and optimization variables used in the stochastic program SP2 for patient assignment, where $\mathcal{U}_{\geq 2}=\{2,\dots, s\}$.

The stochastic program for patient assignment is defined in SP2. Again, unless indicated otherwise, the index variables are assumed to lie in the domains indicated in Table~\ref{table: model definitions3}.
\renewcommand{\lpforall}[1]{&& \forall #1}
\begin{lpformulation}[\text{(SP2)}]
\lplabel{ILP5}
\lpobj[ILP5 obj]{min}{\frac{1}{I}\textstyle\sum_{u=1}^s\sum_{h=1}^H\sum_{i=1}^{I}\sum_{j=0}^J p^{(j)} o_{h,t+u}^{(i,j)}}
\lpeq*{\text{\it First stage:}}{}{}
\lpeq[allocation_today]{\textstyle\sum_{h=1}^H w_{h,t}^{(j)}=1}{j,}
\lpeq[binary domain constraint]{w^{(j)}_{h,t}
\in\{0,1\}}{h, j,}\\\nonumber
\lpeq*{\text{\it Second stage:}}{}{}
\lpeq[allocation WAS]{\textstyle\sum_{h=1}^H x_{h,t+u}^{(i,j)}=A_{t+u}^{(i)}}{u\in\mathcal{U}_{\geq 2},i,j, }
\lpeq[determine occupancy 2]{N_{h,t+u}^{(i,j)}=\tilde{N}_{h,t+u}^{(i)}+\textstyle\sum_{\ell=2}^{u} x_{h,t+\ell}^{(i,j)}+\sum_{k=1}^j w_{h,t}^{(k)}}{u\in\mathcal{U}_+,h,i,j,}
\lpeq[determine overbeds 2]{N^{(i,j)}_{h,t+u}-C_{h,t+u-1}\leq o_{h,t+u}^{(i,j)}}{u\in\mathcal{U}_+,h,i,j,}
\lpeq[integer domain constraints]{N_{h,t+u}^{(i,j)},\,o_{h,t+u}^{(i,j)}\in\mathbb{N}_0}{u\in\mathcal{U}_+,h,i,j,}
\lpeq[integer domain constraints2]{x_{h,t+u}^{(i,j)}\in\mathbb{N}_0}{u\in\mathcal{U}_{\geq2},h,i,j.}
\end{lpformulation}

In the above, Constraint set~\eqref{allocation_today} ensures that all regionally arriving infectious patients in need of assignment today are assigned to a hospital. Constraint set~\eqref{binary domain constraint} considers domain constraints.
Constraint set~\eqref{allocation WAS} ensures that all regionally arriving infectious patients in need of assignment are assigned to a hospital for all scenarios and future days.
Constraint set~\eqref{determine occupancy 2} determines the total occupancy at a hospital based on the assigned regional patients and the autonomous occupancy under the scenario on a given day. 
Constraint set~\eqref{determine overbeds 2} determines the amount of overbeds at a hospital under a scenario on a given day.
Constraint sets~\eqref{integer domain constraints} and~\eqref{integer domain constraints2} are domain constraints.

\begin{table}[htbp]
\small
\caption{Indices, parameters, and decision variables used in the stochastic program SP2 for patient assignment.}
\renewcommand{\arraystretch}{1.3}
\centering
\begin{tabular}{@{} l l @{}}
\toprule
\textbf{Symbol} & \textbf{Description} \\
\midrule
\addlinespace[1ex]
\textbf{Indices}  & \\
$h\in\{1,\dots, H\}$ & Hospital\\
$u\in\{0,\dots, s\}$ & Days ahead \\
$i\in\{1,\dots,I\}$ & First part of scenario:
 occupancy and arrivals after day~$t$\\
$j\in\{0,\dots,J\}$ & Second part of scenario: number of infectious patients to be assigned on day~$t$\\
\addlinespace[1ex]
\textbf{Parameters} & \\
{\it First stage:}\\
$C_{h,t+u}\in\mathbb{N}$ & Number of infectious beds in hospital $h$ on day~$t+u$, $u\in\mathcal{U}_-$ \\
$p^{(j)}\in[0,1]$ & Probability that $j$ regional infectious patients arrive today\\\\
{\it Second stage:}\\
$\tilde{N}_{h,t+u}^{(i)}\in\mathbb{N}$ & Autonomous  occupancy in hospital $h$, at time $t+u$, $u\in\mathcal{U}_+$, scenario $i$ \\
  $A_{t+u}^{(i)}\in\mathbb{N}$   & Patients in need of assignment in $[t+u-1,t+u)$,$u\in\mathcal{U}_+$, scenario $i$\\
\addlinespace[1ex]
\textbf{Variables} & \\
{\it First stage:}\\
$w^{(j)}_{h,t}\in\{0,1\}$ & Indicator that the $j$-th patient at day $t$ is assigned to hospital $h$\\\\
{\it Second stage:}\\
$x_{h,t+u}^{(i,j)}\in\mathbb{N}_0$ &  Patients assigned to hospital $h$ in $[t+u-1,t+u)$, scenario $(i,j)$, $u\in\mathcal{U}_{\geq 2}$\\
$N_{h,t+u}^{(i,j)}\in\mathbb{N}_0$ & Infectious occupancy in hospital $h$ on day~$t+u$, scenario $(i,j)$, $u\in\mathcal{U}_+$\\
$o_{h,t+u}^{(i,j)}\in\mathbb{N}_0$ & Overbeds in hospital $h$ on day~$t+u$ under scenario $(i,j)$, $u\in\mathcal{U}_+$\\
\bottomrule
\end{tabular}
\label{table: model definitions3}
\end{table}
\FloatBarrier

\subsection{Practical implementation of the two linked stochastic programs}\label{sect:relationship_SPs}

Every day, we run SP2 directly after SP1 without knowledge of the actual number of arriving regional infectious patients that day. This construction is especially suited to be used in a simulation to test and evaluate the policies resulting from the two stochastic programs. As mentioned earlier, in reality, one could instead run a stochastic program every time a new regional infectious patient arrives and decide on the placement for this patient alone. This would allow us to use the information on already-arrived autonomous patients and occupancy until then, which would likely lead to an increase in the quality of the solution. However, it would also substantially increase computation time.

Note that in the second stage of SP1, we decide on the assignment of patients to hospitals for today and the upcoming days. In SP2, the assignment decision for today's patients is moved to the first stage, while the second stage still considers the assignment of patients in the upcoming days. Therefore, the quality of the room decisions made by SP1 is dependent on using SP2 for the patient assignment decisions afterward. However, given any room decisions for the next few days, SP2 could be applied independently of SP1 to provide suggestions for patient-hospital assignments. 

\section{Case study: Collaboration of hospitals in Acute Zorg Euregio during the COVID-19 pandemic} \label{sect:results}

This section evaluates the performance of using SP1 and SP2 to determine the number of rooms to open, reserve to open, and assign patients. The performance of this decision rule, denoted~SP, is compared to the performance of two heuristic decision rules in a simulation study.
The decision rules are evaluated based on several \emph{key performance indicators}~(KPIs), namely, the number of overbeds, underbeds, regular care beds used, regular care beds reserved, and the number of rooms added/removed on average per day during the simulation period.
The section ends with a sensitivity analysis. 

The simulation study considers occupancy and regional arrivals by COVID-19 patients in the period 1 October 2021 until 31 December 2021 (91 days) at three hospitals in ROAZ region Acute~Zorg~Euregio~(Euregio): \emph{Medisch~Spectrum~Twente}~(MST), \emph{Ziekenhuisgroep~Twente} (ZGT) and \emph{Streekziekenhuis~Koningin Beatrix}~(SKB).
During COVID-19 outbreaks in the Netherlands, the ROAZ consultation bodies controlled the assignment of infectious patients to hospitals within a region.
At each decision epoch, the method described in Section~\ref{sect:ScenGen} generates scenarios of daily occupancy and arrivals in the period~$[t,t+s)$, used by the considered decision rules. As in~\citet{dijkstra2023dynamic}, the decision epochs are set to 10~AM for each day, which equals the time that hospitals report their occupancy to the region.
At the start of each simulation run, there are no COVID-19 patients at the hospitals. The standard capacity and the capacity of regular care rooms that can be opened for COVID-19 care per hospital can be found in Table~\ref{tabl:capacity}. 
In the following paragraphs, we will first describe each parameter used in the simulation study, after which we will describe the (common) estimation method for this parameter used by the decision rules. 

\begin{table}[htbp]
\centering
\small
\caption{Standard (COVID-19) capacity and regular care room capacity per hospital.}
\label{tabl:capacity}
\begin{tabular}{ccccccl}
\hline
Hospital $h$ & Standard capacity $c_h$ & \multicolumn{5}{l}{Reg. capacity $b_{h,n}$} \\ \hline
MST      & 20                      & 4      & 12      & 6      & 2     & 4   \\
ZGT      & 8                       & 8      & 5       & 5      & 6     & -     \\
SKB      & 8                       & 5      & 7       & -      & -     & -     \\ \hline
\end{tabular}
\end{table}
\FloatBarrier

\begin{figure}[htbp]
\centering
\includegraphics[width=.7\textwidth]{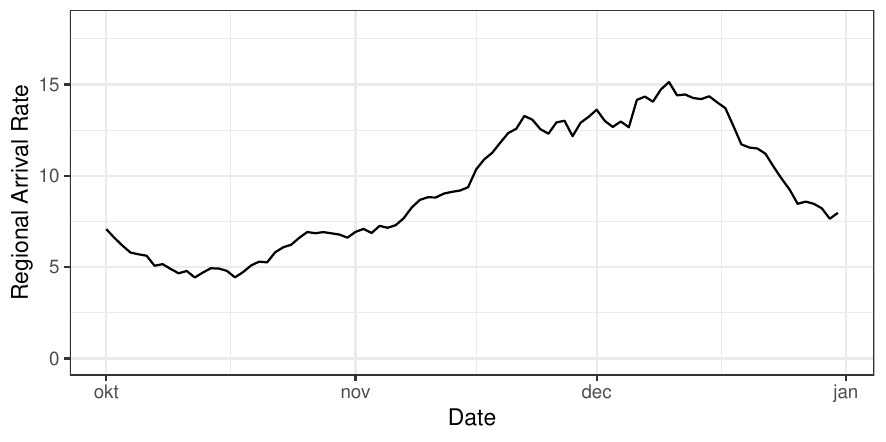}
\caption{Daily regional arrival rates (amount of patients per day) of regional COVID-19 patients in the region for the simulation study of Section~\ref{sect:results} (1/10/2021 - 31/12/2021).}\label{arrival_rates_simstudy_hospitals}
\end{figure}

The input parameters for the simulation study was based on data from the hospital data warehouses, see, e.g.~\citet{Baas2021} for an example. 
Following~\citet{dijkstra2023dynamic}, an exponentially weighted moving average (weight~$0.1$) of the sum of the historic daily COVID-19 arrivals to the three hospitals during the considered period is used to determine the Poisson arrival rate $\lambda_\tau$ of regional infected patients, which is used to sample regional arrivals of COVID-19 patients in the simulation study.  Figure~\ref{arrival_rates_simstudy_hospitals} shows the resulting arrival intensity. 
At each decision epoch $t$, the decision rules use a 5-parameter Richards' curve predictor as described in \citet{Baas2021} and improved in~\citet{dijkstra2023dynamic} to determine $\hat{\lambda}_\tau$ for~$\tau\in[t,t+s)$. 
 The choice was made to use the Richards' curve predictor, and not the daily infections predictor proposed in \citet{dijkstra2023dynamic} as there is no guarantee that in future pandemics, carefully recorded daily infection data will always be publicly available. To stabilize arrival rate predictions at the start of the simulation study, the prediction method for the regional arrival rate uses 2 months of historical daily arrival data prior to the start of the simulation study.
\FloatBarrier

The fractions $f_h$ of COVID-19 patients arriving autonomously at the hospitals are not available from the data and are set to  15\%, 4\%, and 4\% for MST, ZGT, and SKB, respectively. These 
values were determined by pre-computation as they yielded a similar probability of overbeds (around $1\%)$ for all hospitals, given that the hospitals do not scale up (i.e., when only the standard capacity, given in Table~\ref{tabl:capacity}, is utilized).
At each decision epoch~$t$, letting~$t=0$ denote the start of the simulation study, the decision rules use an estimate~$\hat{f}_{h,t}$ of~$f_h$, determined using the formula \begin{equation}
    \hat{f}_{h,t} = (\hat{f}_{h,0} +\textstyle\sum_{u=1}^{t}\tilde{A}_{h,u})/(1+\sum_{u=1}^{t} [A_{u} + \sum_{h=1}^H\tilde{A}_{h,u}] ),\label{def_est_f}
\end{equation}
where $\tilde{A}_{h,u} $ is the number of patients arriving autonomously to hospital $h$ in $[u-1,u)$ of the simulation and $\hat{f}_{h,0}$ is a prior fraction to stabilize the estimation procedure at the start of the simulation period, set to $20\%$, $5\%$, and $5\%$ for MST, ZGT, and SKB. Note that $\hat{f}_{h,t}$ provides a consistent estimator of $f_{h}$ under the model in Section~\ref{sect:ModOcc}. Hence, although the (arbitrary) choice  $20\%$, $5\%$, and $5\%$ for the prior fractions is different from the actual values 15\%, 4\%, and 4\%, the fractions will, in practice, quickly converge to the truth.

A shared LoS distribution in the ward is determined based on pooled data collected from the hospital data warehouses of ZGT~and~SKB, where the weight of the LoS data from ZGT is set to 2.0 as, after consulting with representatives of MST,  it was determined to have the best correspondence with the LoS distribution at MST during this period, where this data was missing. 
 The LoS distribution given as input to the simulation study was determined using the Kaplan-Meier estimator in order to account for right censoring due to patient transfers to other hospitals and patients still residing at the moment of estimation. 
  The LoS distribution is assumed to be known to the decision rules.


We now describe the decision rules considered in the simulation study:
\begin{itemize}
    \item {\bf Individual Hospitals (IH)}:  
    This decision rule mimics the situation in which all regional patients arrive autonomously to the hospitals and hospitals open rooms according to individual forecasts of the infectious occupancy in the coming days. The fractions of regional patients arriving autonomously at the hospitals are set to  48\%, 32\%, and  20\% for MST, ZGT, and SKB, respectively, in agreement with the ratio between maximum capacity available for infectious patients for these hospitals, as given in Table~\ref{tabl:capacity}.
    
    A forecast of the occupancy for the coming days is generated by taking the maximum of the current occupancy at the hospital and a 90\% quantile over scenarios of the occupancy at the hospital on \hbox{day~$t+2$}, as rooms can only be opened in 2 days. 
    The number of rooms to be opened at the hospital is determined as the smallest number of rooms such that the capacity exceeds this occupancy forecast when these rooms are opened. 
    The number of rooms opened and scheduled to open is set to the minimum of the previous decision made by the decision rule and the currently determined amount of rooms to be opened. 
    In order to dampen the fluctuation in decisions over time, the number of opened rooms cannot be decreased by IH when the occupancy forecast for the hospital plus a safety margin exceeds the current hospital capacity. This safety margin is set to 3 beds for MST and 2 beds for ZGT~and~SKB. 
    
    \item {\bf Pandemic Unit (PU)}: 
    This decision rule mimics the situation in which all patients are first sent to the hospital with the largest capacity in the region, which was MST in this case (see Table~\ref{tabl:capacity}).
 When assigning patients, PU sends all patients to MST until the hospital is at capacity minus a safety margin.
    The safety margin equals the occupancy forecast (maximum of occupancy today and $90\%$ quantile of occupancy in 2 days) coming from autonomous arrivals. The remainder of the patients to be assigned are assigned at random with probabilities~61\% vs. 39\% for ZGT~and~SKB, respectively, in agreement with the ratio between maximum capacity available for infectious patients for these hospitals, as given in Table~\ref{tabl:capacity}. 
    
    When determining which rooms to open or close, all arrivals are assumed to go to MST  when this hospital has not opened all its rooms. Otherwise,  the predicted arrival rate of patients to be assigned to MST is set such that MST has no capacity left when the expected amount of arriving patients stay at MST during the interval~$[t,t+s).$
    The occupancy forecast, determined in the same manner as for IH, is included in the calculation of the remaining capacity. It is assumed that the remainder of the regional arrival rate of patients in need of assignment is divided over ZGT~and~SKB in the ratio 61\% vs. 39\%. As in IH, the number of rooms to be opened at each hospital is determined as the smallest number of rooms such that the capacity exceeds the corresponding occupancy forecast when these rooms are opened. The number of rooms opened and scheduled to open is set to the minimum of the previous choice and the currently determined amount of rooms to be opened.  For PU, no rooms at ZGT~and~SKB can be opened or scheduled to be open if MST has not opened all its rooms.  The safety margins for scaling down, used in the same manner as for IH, are set to 3 beds for MST and 2 beds for ZGT~and~SKB.

    \item {\bf Stochastic Program (SP)}: 
    Three cost settings for SP are considered, where the cost vector corresponding to the first, second, and third row of Table~\ref{table_scen_SP} is denoted by \mbox{SP-O,} SP-B, and SP-R, respectively.
    First, the cost parameters of SP-O are chosen such that it outperforms PU (in terms of all cost components given in~\eqref{ILP_shorter obj}) and has a low amount of regional overbeds. 
   Second, the cost parameters of SP-B are chosen such that it outperforms IH and allocates a low amount of regular beds to infectious care. Third, the cost parameters of SP-R are chosen such that it outperforms PU and opens/closes a low amount of regular rooms for infectious care.
    As we will mainly focus on opening or closing a room in the evaluation, we set $\alpha=\beta$ and furthermore set $\gamma=1$ to fix the scale. The parameter choice is found by hand and was completed when changing any of the parameters significantly either led to the regional overbeds increasing or worse performance in the other KPIs than the compared decision rule (IH or PU) in any of the KPIs. 
    In practice, the cost parameters can be based on expert opinion, reflect actual real-life costs, or can be tuned based on scenario analyses.
  \end{itemize}
\begin{table}[htbp]
\centering
\small
\caption{Scenarios for cost vectors considered for SP.}\label{table_scen_SP}
\begin{tabular}{lccccc}
\hline
Setting & $\alpha$ & $\beta$ & $\gamma$ & $\delta$ & $\epsilon$ \\ \hline
Overbeds SP (SP-O)       & 15       & 15      & 1        & 1.5        & 40         \\
Reg. beds SP (SP-B)      & 6       & 6      & 1        & 1        & 13        \\
Open/close rooms (SP-R)     & 60        & 60       & 1        & 1.25        & 25        \\\hline
\end{tabular}
\end{table}
\FloatBarrier
For each decision rule, 250~independent simulation runs of 91~days are performed, and 100~scenarios of regional arrivals and daily occupancy are used by the decision rules for each day in a simulation run. On average, a simulation run for the IH~and~PU decision rules took around 550 seconds, while a simulation run with SP took around 700-1300 seconds.

\subsection{Numerical comparison of decision rules}\label{sect:numerical_comparison_DR}

Table~\ref{tabl:main_result_table} presents the main results of the simulation study.   The underbeds KPI equals the number of regular care beds opened for infectious care that is unoccupied on a given day, averaged over days in a simulation run. This KPI is not directly minimized in the objective function of SP1. To quantify the significance of the differences seen in the table, 95\% \emph{confidence intervals}~(CIs) for the means, based on Student's t-distribution, are also shown in the table. 

The first column of Table~\ref{tabl:main_result_table}, corresponding to IH, shows 
a similar number of overbeds for MST, ZGT, and~SKB, while the number of underbeds is slightly higher for MST and ZGT than for SKB. The ordering in the number of beds reserved and rooms added/removed roughly follows the ordering in the total capacities for the hospitals, while ZGT has the highest need for extra capacity under IH, with around 14 extra beds used on average per day. 

In comparison to IH, the second column of Table~\ref{tabl:main_result_table}, corresponding to PU, shows that the match between occupancy and capacity is worse for MST, as can be seen from the higher amount of underbeds per day, while this match is better for ZGT~and~SKB.  For ZGT~and~SKB, most KPIs decrease, often significantly (as the CIs do not overlap), when comparing IH~to~PU. For MST, the number of beds reserved and rooms added/removed decreases. Looking at the regional level, it is seen that the match between occupancy and capacity becomes better than under IH, as both the average number of over and underbeds decrease at the cost of a higher number of regular care beds used over time. In agreement with the KPIs for MST, the number of beds reserved and rooms added/removed also decreases, indicating less fluctuation in the decisions for scaling up/down COVID-19 capacity over time.
The reason for this can be that under PU, there is less variability in the occupancy at ZGT~and~SKB, as most COVID-19 patients will first go to MST. Due to this, the forecast quality of future occupancy improves and, with it, the performance of PU. 

   The third column of Table~\ref{tabl:main_result_table}, corresponding to SP-O, shows that in comparison to PU, the number of overbeds at ZGT~and~SKB is lower, while the amount of overbeds for MST is higher. The number of regular care beds used and rooms added/removed is more similar for MST and ZGT in comparison to PU. SP-O outperforms IH  for all KPIs and hospitals. The ordering seen in the over and underbeds corresponds to the order of the sizes of the hospitals. 
The similarity of the amount of opened/closed beds for MST and ZGT can be explained by the similarity in total regular care capacity available for COVID-19 for these hospitals, 28 in comparison to 24 respectively (see Table~\ref{tabl:capacity}). 

The bottom rows of the third column Table~\ref{tabl:main_result_table} show the simulated average daily cost for IH,~PU,~and SP-O under the SP-O cost vector, as well as the forecast cost for SP-O (1,~2,~and 3~
days ahead) using the scenarios generated under~SP-O.
The daily costs are realizations of the cost during one simulation run, while the forecasts are made by averaging several scenarios used by SP-O during the simulation run. 
It is seen that IH has the highest cost, then PU and the smallest cost is attained for SP-O, which could be a result of PU and SP-O being geared towards the reduction of regional overbeds, while IH is geared towards the reduction of overbeds at each hospital separately. 
It is seen that the forecast cost increases with the forecast horizon, overshooting the realized cost around~7 for three days ahead. This could be due to the assumption that patients assigned to the hospitals stay there 
 during the period~$[t,t+s)$, while in practice, the LoS might be shorter. 

 The fourth column of Table~\ref{tabl:main_result_table}, corresponding to SP-B, shows that in comparison to IH, the amount of overbeds is higher at MST, while this amount is lower at ZGT~and~SKB. Most other KPIs significantly decrease in comparison to IH. In comparison to PU, the most notable change is the number of regular beds used at MST, which decreases by roughly~$50\%$, while the number of regular beds used at the other hospitals increases. The regional KPIs show that the amount of regional overbeds is close to that of IH, while all other KPIs are significantly lower. In comparison to PU, the number of underbeds and regular beds used is lower. As the regular beds used are the main component of the cost and this KPI involves less uncertainty, the cost forecast stays closer to the realization for SP-B than for SP-O. 

The fifth column of Table~\ref{tabl:main_result_table}, corresponding to SP-R, shows that in comparison to PU, the number of overbeds increases at MST while all other KPIs decrease. For ZGT~and~SKB, the number of underbeds and regular beds used per day increases in comparison to PU, while overbeds, reserved beds, and opened/closed rooms decrease.  Regionally, SP-R outperforms PU and, hence also, IH in all KPIs.  As the cost given to overbeds is higher for SP-R than for SP-B, it is again seen that the forecast becomes higher than the realization for longer horizons.

The results show that SP is a flexible decision rule, where the cost parameters can be tuned to make a trade-off between relevant KPIs, while SP still outperforms rule-of-thumb heuristics. 

\FloatBarrier
\begin{table}[h!]\centering\small
\caption{Average KPIs for each hospital and the total region, along with 95\%~confidence~interval, for the decision rules compared in the simulation study. The average costs for IH~and~PU shown in columns 4, 5, and~6 (rows 21-23) are determined using the same cost parameters as those used for the SP decision rules. }\label{tabl:main_result_table}
\begin{adjustbox}{angle=90}

\begin{tabular}{llccccc}
\hline
Hospital & KPI & IH & PU & SP-O & SP-B & \multicolumn{1}{c}{SP-R} \\ 
\hline
\nopagebreak MST & \nopagebreak Overbeds  & $\text{\phantom{0}0.230 $\pm$ 0.025}$ & $\text{\phantom{0}0.102 $\pm$ 0.012}$ & $\text{\phantom{0}0.140 $\pm$ 0.018}$ & $\text{\phantom{0}0.249 $\pm$ 0.022}$ & $\text{\phantom{0}0.167 $\pm$ 0.019}$ \\
 & \nopagebreak Underbeds  & $\text{\phantom{0}4.654 $\pm$ 0.135}$ & $\text{\phantom{0}6.570 $\pm$ 0.073}$ & $\text{\phantom{0}3.729 $\pm$ 0.104}$ & $\text{\phantom{0}2.747 $\pm$ 0.083}$ & $\text{\phantom{0}3.586 $\pm$ 0.110}$ \\
 & \nopagebreak Reg. beds used  & $\text{13.056 $\pm$ 0.178}$ & $\text{22.903 $\pm$ 0.115}$ & $\text{12.187 $\pm$ 0.188}$ & $\text{10.946 $\pm$ 0.179}$ & $\text{11.588 $\pm$ 0.206}$ \\
 & \nopagebreak Beds reserved  & $\text{\phantom{0}1.163 $\pm$ 0.034}$ & $\text{\phantom{0}0.835 $\pm$ 0.024}$ & $\text{\phantom{0}0.828 $\pm$ 0.024}$ & $\text{\phantom{0}0.926 $\pm$ 0.026}$ & $\text{\phantom{0}0.711 $\pm$ 0.022}$ \\
 & \nopagebreak Rooms added/removed  & $\text{\phantom{0}0.119 $\pm$ 0.004}$ & $\text{\phantom{0}0.078 $\pm$ 0.003}$ & $\text{\phantom{0}0.085 $\pm$ 0.003}$ & $\text{\phantom{0}0.100 $\pm$ 0.003}$ & $\text{\phantom{0}0.057 $\pm$ 0.001}$ \\
\rule{0pt}{1.7\normalbaselineskip}ZGT & \nopagebreak Overbeds  & $\text{\phantom{0}0.269 $\pm$ 0.033}$ & $\text{\phantom{0}0.181 $\pm$ 0.025}$ & $\text{\phantom{0}0.096 $\pm$ 0.012}$ & $\text{\phantom{0}0.165 $\pm$ 0.013}$ & $\text{\phantom{0}0.106 $\pm$ 0.011}$ \\
 & \nopagebreak Underbeds  & $\text{\phantom{0}4.779 $\pm$ 0.093}$ & $\text{\phantom{0}3.359 $\pm$ 0.084}$ & $\text{\phantom{0}3.397 $\pm$ 0.079}$ & $\text{\phantom{0}2.573 $\pm$ 0.060}$ & $\text{\phantom{0}3.660 $\pm$ 0.077}$ \\
 & \nopagebreak Reg. beds used  & $\text{14.309 $\pm$ 0.128}$ & $\text{\phantom{0}9.126 $\pm$ 0.118}$ & $\text{12.310 $\pm$ 0.171}$ & $\text{11.141 $\pm$ 0.168}$ & $\text{12.877 $\pm$ 0.164}$ \\
 & \nopagebreak Beds reserved  & $\text{\phantom{0}1.001 $\pm$ 0.029}$ & $\text{\phantom{0}0.879 $\pm$ 0.026}$ & $\text{\phantom{0}0.846 $\pm$ 0.024}$ & $\text{\phantom{0}0.906 $\pm$ 0.026}$ & $\text{\phantom{0}0.779 $\pm$ 0.022}$ \\
 & \nopagebreak Rooms added/removed  & $\text{\phantom{0}0.085 $\pm$ 0.003}$ & $\text{\phantom{0}0.089 $\pm$ 0.003}$ & $\text{\phantom{0}0.084 $\pm$ 0.003}$ & $\text{\phantom{0}0.088 $\pm$ 0.003}$ & $\text{\phantom{0}0.047 $\pm$ 0.001}$ \\
\rule{0pt}{1.7\normalbaselineskip}SKB & \nopagebreak Overbeds  & $\text{\phantom{0}0.207 $\pm$ 0.029}$ & $\text{\phantom{0}0.268 $\pm$ 0.036}$ & $\text{\phantom{0}0.090 $\pm$ 0.010}$ & $\text{\phantom{0}0.159 $\pm$ 0.014}$ & $\text{\phantom{0}0.110 $\pm$ 0.012}$ \\
 & \nopagebreak Underbeds  & $\text{\phantom{0}3.388 $\pm$ 0.107}$ & $\text{\phantom{0}1.986 $\pm$ 0.078}$ & $\text{\phantom{0}2.509 $\pm$ 0.071}$ & $\text{\phantom{0}1.906 $\pm$ 0.053}$ & $\text{\phantom{0}2.426 $\pm$ 0.068}$ \\
 & \nopagebreak Reg. beds used  & $\text{\phantom{0}7.174 $\pm$ 0.101}$ & $\text{\phantom{0}4.836 $\pm$ 0.069}$ & $\text{\phantom{0}7.221 $\pm$ 0.144}$ & $\text{\phantom{0}6.816 $\pm$ 0.130}$ & $\text{\phantom{0}7.381 $\pm$ 0.135}$ \\
 & \nopagebreak Beds reserved  & $\text{\phantom{0}0.679 $\pm$ 0.023}$ & $\text{\phantom{0}0.437 $\pm$ 0.019}$ & $\text{\phantom{0}0.469 $\pm$ 0.018}$ & $\text{\phantom{0}0.486 $\pm$ 0.019}$ & $\text{\phantom{0}0.417 $\pm$ 0.016}$ \\
 & \nopagebreak Rooms added/removed  & $\text{\phantom{0}0.060 $\pm$ 0.003}$ & $\text{\phantom{0}0.040 $\pm$ 0.002}$ & $\text{\phantom{0}0.042 $\pm$ 0.002}$ & $\text{\phantom{0}0.042 $\pm$ 0.002}$ & $\text{\phantom{0}0.032 $\pm$ 0.001}$ \\
\rule{0pt}{1.7\normalbaselineskip}Region & \nopagebreak Overbeds  & $\text{\phantom{0}0.706 $\pm$ 0.051}$ & $\text{\phantom{0}0.552 $\pm$ 0.046}$ & $\text{\phantom{0}0.326 $\pm$ 0.035}$ & $\text{\phantom{0}0.573 $\pm$ 0.041}$ & $\text{\phantom{0}0.383 $\pm$ 0.037}$ \\
 & \nopagebreak Underbeds  & $\text{12.821 $\pm$ 0.197}$ & $\text{11.915 $\pm$ 0.137}$ & $\text{\phantom{0}9.635 $\pm$ 0.181}$ & $\text{\phantom{0}7.226 $\pm$ 0.146}$ & $\text{\phantom{0}9.672 $\pm$ 0.193}$ \\
 & \nopagebreak Reg. beds used  & $\text{34.539 $\pm$ 0.251}$ & $\text{36.865 $\pm$ 0.208}$ & $\text{31.717 $\pm$ 0.240}$ & $\text{28.902 $\pm$ 0.242}$ & $\text{31.846 $\pm$ 0.233}$ \\
 & \nopagebreak Beds reserved  & $\text{\phantom{0}2.843 $\pm$ 0.052}$ & $\text{\phantom{0}2.151 $\pm$ 0.039}$ & $\text{\phantom{0}2.143 $\pm$ 0.042}$ & $\text{\phantom{0}2.318 $\pm$ 0.043}$ & $\text{\phantom{0}1.907 $\pm$ 0.034}$ \\
 & \nopagebreak Rooms added/removed  & $\text{\phantom{0}0.263 $\pm$ 0.006}$ & $\text{\phantom{0}0.207 $\pm$ 0.005}$ & $\text{\phantom{0}0.211 $\pm$ 0.004}$ & $\text{\phantom{0}0.230 $\pm$ 0.005}$ & $\text{\phantom{0}0.136 $\pm$ 0.001}$ \\
 & \rule{0pt}{1.7\normalbaselineskip}Avg. cost IH  & $\text{\phantom{0.0000}-\phantom{00.00}}$ & $\text{\phantom{0.0000}-\phantom{00.00}}$ & $\text{70.990 $\pm$ 2.077}$ & $\text{48.139 $\pm$ 0.726}$ & $\text{71.547 $\pm$ 1.329}$ \\
 & \nopagebreak Avg. cost PU  & $\text{\phantom{0.0000}-\phantom{00.00}}$ & $\text{\phantom{0.0000}-\phantom{00.00}}$ & $\text{65.291 $\pm$ 1.870}$ & $\text{47.439 $\pm$ 0.654}$ & $\text{65.797 $\pm$ 1.202}$ \\
 & \nopagebreak Avg. cost SP  & $\text{\phantom{0.0000}-\phantom{00.00}}$ & $\text{\phantom{0.0000}-\phantom{00.00}}$ & $\text{51.129 $\pm$ 1.407}$ & $\text{40.055 $\pm$ 0.598}$ & $\text{51.959 $\pm$ 0.932}$ \\
 & \nopagebreak Forecast cost SP (hor. 1)  & $\text{\phantom{0.0000}-\phantom{00.00}}$ & $\text{\phantom{0.0000}-\phantom{00.00}}$ & $\text{48.990 $\pm$ 1.238}$ & $\text{38.469 $\pm$ 0.508}$ & $\text{50.248 $\pm$ 0.784}$ \\
 & \nopagebreak Forecast cost SP (hor. 2)  & $\text{\phantom{0.0000}-\phantom{00.00}}$ & $\text{\phantom{0.0000}-\phantom{00.00}}$ & $\text{50.908 $\pm$ 1.325}$ & $\text{39.255 $\pm$ 0.501}$ & $\text{51.855 $\pm$ 0.831}$ \\
 & \nopagebreak Forecast cost SP (hor. 3)  & $\text{\phantom{0.0000}-\phantom{00.00}}$ & $\text{\phantom{0.0000}-\phantom{00.00}}$ & $\text{58.380 $\pm$ 1.797}$ & $\text{41.768 $\pm$ 0.630}$ & $\text{59.388 $\pm$ 1.127}$ \\
\hline 
\end{tabular}

\end{adjustbox}
\end{table}
\FloatBarrier

\subsection{Sensitivity analyses}\label{sect:sens_analysis}

In this section, we analyze the sensitivity of the decision rules IH, PU, and SP to several parameters. In order to analyze SP, we chose to analyze only SP-O  as the quality of the solution depends most on the scenarios given as input to SP.  First, we analyze the performance of the IH~and~PU decision rules when changing the quantile used for determining the occupancy forecast to investigate the robustness of the cost vectors in Table~\ref{table_scen_SP} to the chosen quantile.
Second, we analyze the sensitivity of SP-O to the lookahead horizon, the number of scenarios, as well as the arrival rate predictor used in generating the scenarios. Third, we compare SP-O with deterministic versions of SP-O based on the median and a quantile over scenarios to investigate whether results similar to those for SP can be reached with a less computationally intensive approach.

\subsubsection{Sensitivity of IH~and~PU to choice of quantile}

Table~\ref{tab:comparison_IHquantiles} shows the results for IH~and~PU when the $80\%,\,85\%,\,90\%$, and $95\%$ quantile is taken over scenarios to determine the forecasts. The number of overbeds decreases for IH when the quantile over scenarios increases. For IH, the number of underbeds and regular beds used increases, while the number of beds reserved and rooms added/removed slightly decreases. PU seems less sensitive to the choice of quantile, with most differences not being significant. The main difference between IH~and~PU is that the number of beds reserved increases slightly for PU in the value of the quantile, while it decreases for IH.  
Tables~\ref{tabl:main_result_table}~and~\ref{tab:comparison_IHquantiles} show that all SP decision rules still outperform the respective heuristic~(IH or PU) even though they were designed to outperform it for other values of the quantile. 
\FloatBarrier
\begin{table}[htbp]\centering
\small
\caption{Regional KPIs for decision rules IH~and~PU using the $80,\,85,\,90,\,95\%$ quantile over scenarios to forecast the occupancy in 2 days. }\hspace{-3mm}\begin{tabular}{llcccc}
\hline
DR & KPI & 80\% & 85\% & 90\% & \multicolumn{1}{c}{95\%} \\ 
\hline
\nopagebreak IH & \nopagebreak Overbeds  & $\text{\phantom{0}0.971 $\pm$ 0.053}$ & $\text{\phantom{0}0.831 $\pm$ 0.052}$ & $\text{\phantom{0}0.706 $\pm$ 0.051}$ & $\text{\phantom{0}0.575 $\pm$ 0.051}$ \\
 & \nopagebreak Underbeds  & $\text{10.398 $\pm$ 0.168}$ & $\text{11.489 $\pm$ 0.178}$ & $\text{12.821 $\pm$ 0.197}$ & $\text{14.690 $\pm$ 0.220}$ \\
 & \nopagebreak Reg. beds used  & $\text{31.954 $\pm$ 0.247}$ & $\text{33.134 $\pm$ 0.248}$ & $\text{34.539 $\pm$ 0.251}$ & $\text{36.476 $\pm$ 0.240}$ \\
 & \nopagebreak Beds reserved  & $\text{\phantom{0}2.908 $\pm$ 0.052}$ & $\text{\phantom{0}2.875 $\pm$ 0.054}$ & $\text{\phantom{0}2.843 $\pm$ 0.052}$ & $\text{\phantom{0}2.841 $\pm$ 0.054}$ \\
 & \nopagebreak Rooms add./rem.  & $\text{\phantom{0}0.272 $\pm$ 0.006}$ & $\text{\phantom{0}0.267 $\pm$ 0.007}$ & $\text{\phantom{0}0.263 $\pm$ 0.006}$ & $\text{\phantom{0}0.256 $\pm$ 0.006}$ \\
\rule{0pt}{1.7\normalbaselineskip}PU & \nopagebreak Overbeds  & $\text{\phantom{0}0.564 $\pm$ 0.045}$ & $\text{\phantom{0}0.559 $\pm$ 0.045}$ & $\text{\phantom{0}0.552 $\pm$ 0.046}$ & $\text{\phantom{0}0.526 $\pm$ 0.044}$ \\
 & \nopagebreak Underbeds  & $\text{11.802 $\pm$ 0.132}$ & $\text{11.841 $\pm$ 0.134}$ & $\text{11.915 $\pm$ 0.137}$ & $\text{12.126 $\pm$ 0.139}$ \\
 & \nopagebreak Reg. beds used  & $\text{36.743 $\pm$ 0.207}$ & $\text{36.784 $\pm$ 0.207}$ & $\text{36.865 $\pm$ 0.208}$ & $\text{37.091 $\pm$ 0.208}$ \\
 & \nopagebreak Beds reserved  & $\text{\phantom{0}2.128 $\pm$ 0.039}$ & $\text{\phantom{0}2.140 $\pm$ 0.039}$ & $\text{\phantom{0}2.151 $\pm$ 0.039}$ & $\text{\phantom{0}2.158 $\pm$ 0.039}$ \\
 & \nopagebreak Rooms add./rem.  & $\text{\phantom{0}0.206 $\pm$ 0.005}$ & $\text{\phantom{0}0.207 $\pm$ 0.005}$ & $\text{\phantom{0}0.207 $\pm$ 0.005}$ & $\text{\phantom{0}0.204 $\pm$ 0.005}$ \\
\hline 
\end{tabular}

\label{tab:comparison_IHquantiles}
\end{table}
\FloatBarrier
\subsubsection{Comparison of lookahead horizons, number of scenarios, and arrival predictors}
Table~\ref{tab:comparison_horizon} shows the KPIs for SP-O, as well as the resulting KPIs for using SP with a different arrival rate predictor and different lookahead horizons.
The second column shows the KPIs when taking the $90\%$ upper bound of the confidence interval for the predicted arrival rate, following from ordinary least squares, to generate scenarios of regional arrivals and occupancy in the stochastic program. For this arrival rate predictor, the table shows that SP-O results in a lower amount of overbeds, beds reserved, and rooms added/removed, while it results in a higher amount of regular beds used, which makes sense as a higher arrival rate leads to a higher forecast occupancy and required capacity. 
Tables~\ref{tabl:main_result_table} and~\ref{tab:comparison_horizon} show that SP-O with the UB arrival rate predictor still outperforms PU, while the average cost for SP-O under the UB arrival rate predictor is lower, although not significantly lower based on the CIs, indicating that this arrival rate could be used as an alternative to the arrival rate used in Table~\ref{tabl:main_result_table}. The lower cost under the UB arrival rate could possibly be explained by the negative bias for the Richards' curve predictor observed in~\citet{Baas2021}~and~\citet{dijkstra2023dynamic}.  
The table shows that the error between realized and forecast cost grows with the horizon of the scenarios, which is expected, as the predictor used for generating arrival scenarios lies above that used for SP-O.

\begin{figure}[htp]
	\centering
	\includegraphics[width = 0.9\textwidth]{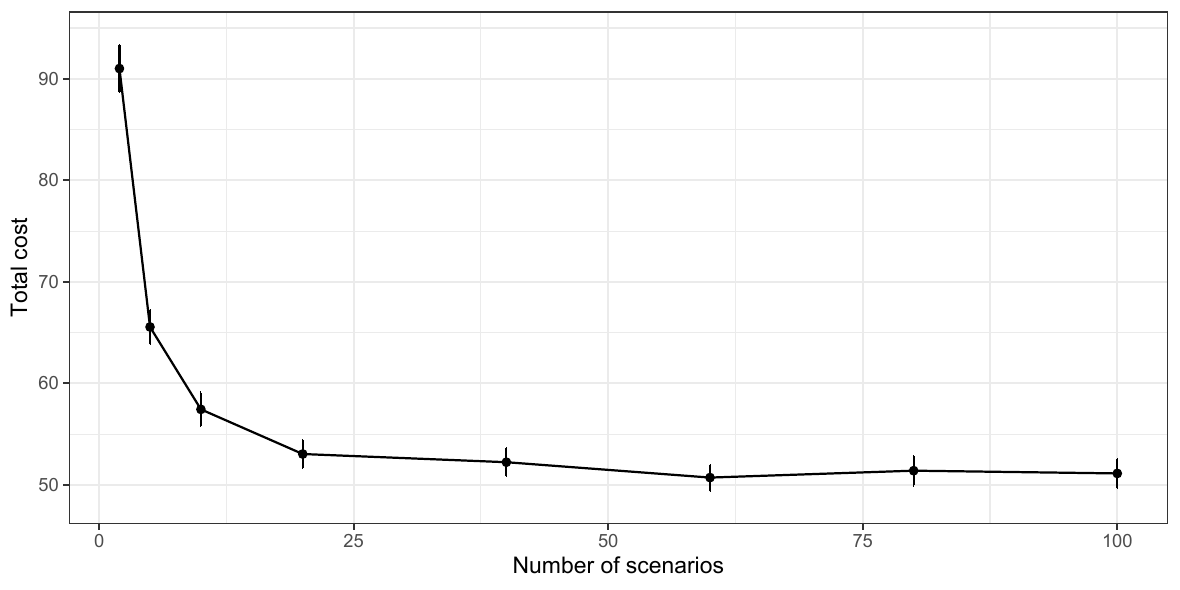}
	\caption{Average total cost (and 95\% CI) for SP-O vs. the number of included scenarios.}\label{plot_cost_vs_nscen}
\end{figure}

The last two columns of Table~\ref{tab:comparison_horizon} show results for different lookahead horizons $s=4,3$. 
The table shows that, on average, the realized daily average amount of overbeds increases, and the realized average daily amount of underbeds and regular beds used decreases when decreasing the horizon for the scenarios.
For shorter scenario horizons, the stochastic program often underestimates the total cost, where for a scenario horizon of 3 days, the forecast cost stays below the realized cost for all forecast horizons. 
As the difference in average cost between SP-O  with a horizon of 4 or 5~days is not significant, it is concluded that taking a scenario horizon of 4~days would also have been sufficient. Comparing the average costs in columns 1 and 2 and columns 1, 3, and 4 in Table~\ref{tab:comparison_horizon}, it is concluded that the program is more sensitive to the scenario horizon than the arrival rate predictor used. 
\FloatBarrier
     \begin{table}[htbp]
\centering\small
\caption{Regional KPIs for SP-O for lookahead horizons 5~(column SP-O), 4, and 3 days and when taking the $90\%$ CI upper bound for the predicted arrival intensity (UB arrival rate). }\hspace{-3mm}\begin{tabular}{lcccc}
\hline
KPI  & SP-O & UB arrival rate & Horizon 4 days & \multicolumn{1}{c}{Horizon 3 days} \\ 
\hline
\nopagebreak Overbeds  & $\text{\phantom{0}0.326 $\pm$ 0.035}$ & $\text{\phantom{0}0.253 $\pm$ 0.036}$ & $\text{\phantom{0}0.380 $\pm$ 0.038}$ & $\text{\phantom{0}0.559 $\pm$ 0.041}$ \\
\nopagebreak Underbeds  & $\text{\phantom{0}9.635 $\pm$ 0.181}$ & $\text{11.212 $\pm$ 0.210}$ & $\text{\phantom{0}8.832 $\pm$ 0.172}$ & $\text{\phantom{0}7.610 $\pm$ 0.148}$ \\
\nopagebreak Reg. beds used  & $\text{31.717 $\pm$ 0.240}$ & $\text{33.464 $\pm$ 0.240}$ & $\text{30.814 $\pm$ 0.239}$ & $\text{29.183 $\pm$ 0.242}$ \\
\nopagebreak Beds reserved  & $\text{\phantom{0}2.143 $\pm$ 0.042}$ & $\text{\phantom{0}2.023 $\pm$ 0.040}$ & $\text{\phantom{0}2.355 $\pm$ 0.043}$ & $\text{\phantom{0}2.198 $\pm$ 0.039}$ \\
\nopagebreak Rooms added/removed  & $\text{\phantom{0}0.211 $\pm$ 0.004}$ & $\text{\phantom{0}0.187 $\pm$ 0.004}$ & $\text{\phantom{0}0.228 $\pm$ 0.005}$ & $\text{\phantom{0}0.218 $\pm$ 0.004}$ \\
\rule{0pt}{1.7\normalbaselineskip}Avg. cost SP  & $\text{51.129 $\pm$ 1.407}$ & $\text{49.429 $\pm$ 1.459}$ & $\text{52.956 $\pm$ 1.533}$ & $\text{58.120 $\pm$ 1.676}$ \\
\nopagebreak Forecast cost SP (hor. 1)  & $\text{48.990 $\pm$ 1.238}$ & $\text{50.002 $\pm$ 1.423}$ & $\text{50.323 $\pm$ 1.293}$ & $\text{52.878 $\pm$ 1.328}$ \\
\nopagebreak Forecast cost SP (hor. 2)  & $\text{50.908 $\pm$ 1.325}$ & $\text{55.907 $\pm$ 1.644}$ & $\text{50.311 $\pm$ 1.356}$ & $\text{51.105 $\pm$ 1.307}$ \\
\nopagebreak Forecast cost SP (hor. 3)  & $\text{58.380 $\pm$ 1.797}$ & $\text{68.829 $\pm$ 2.325}$ & $\text{57.086 $\pm$ 1.867}$ & $\text{56.332 $\pm$ 1.738}$ \\
\hline 
\end{tabular}

\label{tab:comparison_horizon}
\end{table}
\FloatBarrier
Another setting of SP is the number of scenarios used. Figure~\ref{plot_cost_vs_nscen} shows the average cost of SP-O, as well as 95\%~CI for the expected cost vs. the number of scenarios. It is seen that the average cost stabilizes roughly after using 50 scenarios, indicating that taking 50 scenarios would have also been sufficient for SP-O. 

\subsubsection{Comparison to deterministic programs}
Table~\ref{tab:comparison_deterministic_sps} shows the KPIs for SP-O, as well as the resulting KPIs when using deterministic versions of SP-O. The second and third columns of  Table~\ref{tab:comparison_deterministic_sps} show the KPIs for SP-O when taking the median and $85\%$ quantile over scenarios of arrivals and occupancy and using this as the only scenario. The $85\%$ quantile is chosen as it resulted in the smallest average cost when compared to deterministic programs with quantiles~$60\%,70\%,80\%,85\%,90\%,$ and $95\%$. 

The results show that taking the median over scenarios leads to a more optimistic forecast than under SP-O, resulting in a realized average cost that is more than 50\% higher than the~1-day ahead forecast. The long-term forecasts show an even more extreme behavior, where the realized cost is more than three times higher. Regionally, this leads to a higher amount of realized overbeds and rooms added/removed, while there is a lower amount of regular beds used compared to~SP-O.

\begin{table}[htbp]
\centering
\small
\caption{Regional KPIs for SP-O, as well as the deterministic version of the stochastic programs, when taking the median and $85\%$ quantile over the arrivals and occupancy scenarios for each~day. }
\begin{tabular}{lccc}
\hline
KPI  & SP-O & \;\;\;Median Scen. & \multicolumn{1}{c}{\;\;\;Quantile Scen.} \\ 
\hline
\nopagebreak Overbeds  & $\text{\phantom{0}\phantom{0}0.326 $\pm$ 0.035}$ & $\text{\phantom{0}\phantom{0}2.455 $\pm$ 0.061}$ & $\text{\phantom{0}\phantom{0}0.419 $\pm$ 0.047}$ \\
\nopagebreak Underbeds  & $\text{\phantom{0}\phantom{0}9.635 $\pm$ 0.181}$ & $\text{\phantom{0}\phantom{0}4.084 $\pm$ 0.083}$ & $\text{\phantom{0}12.029 $\pm$ 0.199}$ \\
\nopagebreak Reg. beds used  & $\text{\phantom{0}31.717 $\pm$ 0.240}$ & $\text{\phantom{0}23.744 $\pm$ 0.245}$ & $\text{\phantom{0}34.079 $\pm$ 0.240}$ \\
\nopagebreak Beds reserved  & $\text{\phantom{0}\phantom{0}2.143 $\pm$ 0.042}$ & $\text{\phantom{0}\phantom{0}2.209 $\pm$ 0.034}$ & $\text{\phantom{0}\phantom{0}2.340 $\pm$ 0.044}$ \\
\nopagebreak Rooms added/removed  & $\text{\phantom{0}\phantom{0}0.211 $\pm$ 0.004}$ & $\text{\phantom{0}\phantom{0}0.239 $\pm$ 0.004}$ & $\text{\phantom{0}\phantom{0}0.196 $\pm$ 0.005}$ \\
\rule{0pt}{1.7\normalbaselineskip}Avg. cost SP  & $\text{\phantom{0}51.129 $\pm$ 1.407}$ & $\text{128.850 $\pm$ 2.542}$ & $\text{\phantom{0}57.299 $\pm$ 1.883}$ \\
\nopagebreak Forecast cost SP (hor. 1)  & $\text{\phantom{0}48.990 $\pm$ 1.238}$ & $\text{\phantom{0}85.416 $\pm$ 1.812}$ & $\text{\phantom{0}72.173 $\pm$ 2.988}$ \\
\nopagebreak Forecast cost SP (hor. 2)  & $\text{\phantom{0}50.908 $\pm$ 1.325}$ & $\text{\phantom{0}56.003 $\pm$ 1.257}$ & $\text{105.605 $\pm$ 4.429}$ \\
\nopagebreak Forecast cost SP (hor. 3)  & $\text{\phantom{0}58.380 $\pm$ 1.797}$ & $\text{\phantom{0}41.251 $\pm$ 1.445}$ & $\text{156.446 $\pm$ 5.864}$ \\
\hline 
\end{tabular}

\label{tab:comparison_deterministic_sps}
\end{table}
\FloatBarrier

The program that results when taking the $85\%$ quantiles shows another extreme, where the forecast cost lies above the realized cost, being more than three times higher than the realized cost for a forecast horizon of four days. For this program, the amount of rooms added/removed is slightly lower than that under SP-O, while the other KPIs are higher. The realized average cost is around 12\% higher. 

Even though the deterministic decision rule with a well-chosen quantile percentage comes close to the performance of SP-O, this can only be achieved after tuning the quantile parameter based on the outcomes of the simulation study. In the real world, the cost parameters of SP might correspond to actual costs and can be determined without using simulations, while this is not the case for the choice of quantile over scenarios. Hence, for such situations, we likely overestimate the performance of the deterministic decision rule in our results.

Comparing the results shown in Table~\ref{tab:comparison_deterministic_sps} with those shown in Table~\ref{tab:comparison_horizon}, the results for \mbox{SP-O} change more drastically when making the program deterministic, in comparison to when changing certain settings of the stochastic program such as the arrival rate predictor, horizon and number of scenarios.


\section{Discussion~\label{sect:discussion}}
In this paper, we developed a joint regional decision-making approach for a region with collaborating hospitals to allocate hospital bed capacity to regular or infectious care and to assign infectious patients to hospitals in a region during an infectious outbreak, aiming to serve all infectious demand while maintaining regular care. 
The objective of the approach is to minimize the sum of costs for rooms currently being available, costs for making rooms ready, and costs for opening or closing rooms to accommodate infectious care patients, as well as costs for not being able to accommodate infectious care patients.
The presented decision rules result from a stochastic lookahead approach by solving two stochastic programs with scenarios in sequence on each day, where the first program makes decisions on opening and closing hospital rooms for infectious care, while the second one makes decisions on the allocation of new infectious patients to hospitals in the region. 

In our numerical experiments, we evaluated the performance of the developed decision rules considering the use case of COVID-19 in a region with three hospitals in the Netherlands at the end of 2021. The results indicate that our approach considerably outperforms simpler decision rules for non-collaborating hospitals and for having a pandemic unit, showing a clear benefit of using our solution approach and of hospital collaboration in a region during infectious outbreaks. 
The comparison of one of these heuristics, the pandemic unit, with regional collaboration led to the conclusion for the considered region that a strategy involving regional collaboration with dynamic capacity allocation and patient assignment to all hospitals would be preferred over a pandemic unit.
This, in part, led to the conclusion by Acute Zorg Euregio that, in order to continue access to regular care as much as possible, regional collaboration would be preferred over a pandemic unit. This indicates the practical relevance of this paper; see, e.g.,~\citet{artikel_NOS}.
The numerical results also demonstrate that our stochastic lookahead approach is superior to deterministic lookaheads. For our use case, we show that only a few scenarios in the stochastic programs, i.e., approximately 50-100, and a short lookahead horizon, i.e., 4-5 days, are needed, ensuring a fast runtime of the stochastic programs. Further, a computational advantage of our decision rules is that they do not rely on tunable parameters compared to simpler ones. Therefore, we can directly apply our rules instead of having to tune parameters in a simulation upfront. 

The numerical experiments showed that a small number (around 50) of scenarios are sufficient in SP1 and SP2. When a much larger number of scenarios needs to be taken into account, the standard solution method used to generate solutions of SP1 and SP2 becomes computationally intractable. Solution methods that take the structure of the stochastic programs into account, such as the L-shaped method, can be employed to overcome this potential difficulty.

Our generic solution approach allows the application of decision rules during future infectious outbreaks to guide bed capacity and patient assignment. For that, the weights for the objective function must be determined by the decision-makers, and scenarios for the stochastic programs must be generated. To create those input scenarios, any suitable method to predict infectious arrival rates can be applied together with our method or alternative methods to produce bed occupancy scenarios. Further, the stochastic programs can easily be adjusted to model different numbers of days it takes to empty a room. The constraints encoding the order of opening and closing rooms in a hospital can also easily be removed or adjusted.

For future research, solution quality could potentially be increased further by revisiting some of our simplifying assumptions. First, one could investigate including a stochastic number of days to empty a room, reconsider the assumption that assigned patients stay until the end of the lookahead horizon, that the LoS distribution is known, or that the fractions of autonomously arriving patients are time-invariant.
The assumption that patients stay until the end of the lookahead horizon leads to an overestimation of the occupancy by patients arriving during the forecast period. This will lead to an overestimation of the number of overbeds, and hence, when real-life cost parameters are used in the model, the cost due to overbeds is overestimated by the model. In a sense, SP1 and SP2 are more conservative than when the LoS of future allocated patients are taken into account. This can lead to a too large amount of rooms being made available for infectious patients. Increasing the cost of opening or closing rooms mitigates this effect, as is shown in Section~\ref{sect:numerical_comparison_DR}.   
Second, one could also possibly consider multiple patient types. Third, multi-stage stochastic programs can be employed instead of two-stage programs. Fourth, one could solve a stochastic program for each arriving regional infectious patient instead of one single program for all patients of the day at the beginning of the day. This way, one can include the realizations of the autonomous inflow of patients until the arrival of that patient. However, while it might be possible to work with multi-stage programs and stochastic programs for individual patients at decision times in real life, evaluating those programs via simulation could become intractable. 

In this paper, we considered three choices of the weights of the objective function of the stochastic program, chosen such that one of the considered heuristics is outperformed and the stochastic program results in either a low amount of overbeds, regular beds used, or rooms opened and/or closed. These parameters were found by hand, and in future research, it could be interesting to consider optimizing the weight parameters such that the stochastic program satisfies the aforementioned properties. In practice, at the start of a pandemic, 
such an optimization method can optimize the weights based on scenario analysis, or the weights might first reflect real-life costs, while the weight optimization method can force the behavior to have desired properties if, at interim analyses, the previously set weights seem insufficient to do so.  

To support even more planning decisions during infectious outbreaks, \emph{intensive care unit}~(ICU) capacity and patient transfers between the ward and ICU, as well as between hospitals, could be included. We could further model bed capacity that can be added without reducing regular care bed capacity. 
Furthermore, occupancy by high-priority regular care patients who have to be admitted to a hospital (i.e., regular patients arriving autonomously) could also be modeled according to the proposed queueing model and can be taken into account in our proposed decision rule. 
For instance, one could decide to allocate less infectious patients to hospitals where the strain on (future) capacity by high-priority~(i.e., more severely ill) regular patients will be higher. 
Finally, it would be interesting to extend our modeling approach from the regional to the national level and to include the transfer of patients to other hospitals (either inside or outside a region).

\section{Acknowledgments}

The authors thank Manon Bruens (Acute Zorg Euregio), Bart Veltman, Sophie Ligtenstein, and Franka Schellekens~(Rhythm) for their collaboration in obtaining the input data for the case study in this paper. 

This article was written as part of the European HORIZON project RAPIDE.

“Funded by the European Union. Views and opinions expressed are however those of the author(s) only and do not necessarily reflect those of the European Union or the European Health and Digital Executive Agency (HaDEA). Neither the European Union nor the granting authority can be held responsible for them.” 
\begin{figure}[htp]
	\centering
	\includegraphics[width = 0.5\textwidth]{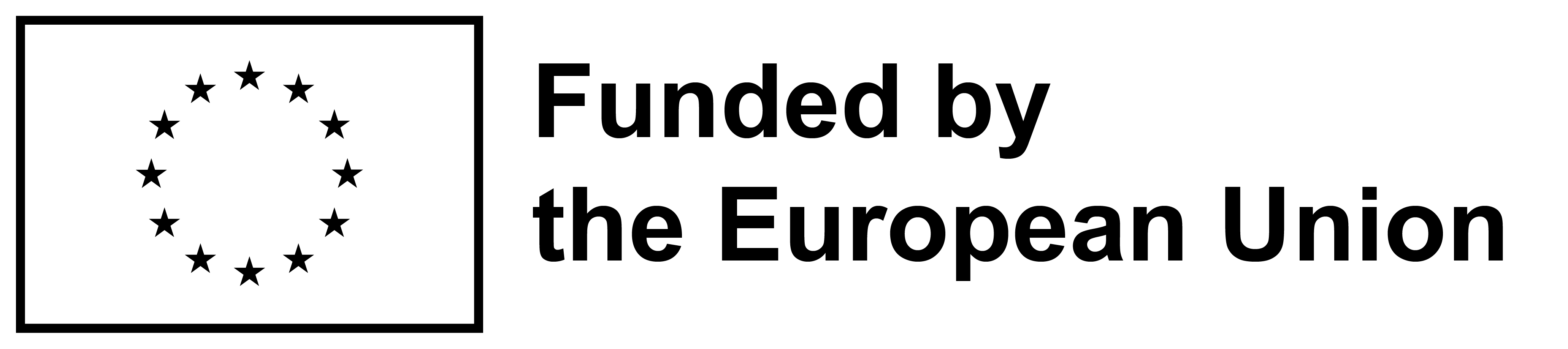}

\end{figure}


\bibliographystyle{elsarticle-harv} 
\bibliography{arxiv_version.bib}
\appendix

\begin{table}
\small
\section{Full tables for Section~5.2}
\caption{ KPIs per hospital for decision rules (decision rule) IH and PU using the $80,\,85,\,90,\,95\%$ quantile over scenarios to forecast the occupancy in 2 days. }
\centering
\begin{adjustbox}{angle=90}
\centering
\small
\begin{tabular}{lllcccc}
\hline
Decision rule & Hospital & KPI & 80\% & 85\% & 90\% & \multicolumn{1}{c}{95\%} \\ 
\hline
\nopagebreak IH & \nopagebreak MST & \nopagebreak Overbeds  & $\text{\phantom{0}0.350 $\pm$ 0.028}$ & $\text{\phantom{0}0.285 $\pm$ 0.026}$ & $\text{\phantom{0}0.230 $\pm$ 0.025}$ & $\text{\phantom{0}0.168 $\pm$ 0.024}$ \\
 &  & \nopagebreak Underbeds  & $\text{\phantom{0}3.641 $\pm$ 0.117}$ & $\text{\phantom{0}4.092 $\pm$ 0.124}$ & $\text{\phantom{0}4.654 $\pm$ 0.135}$ & $\text{\phantom{0}5.439 $\pm$ 0.148}$ \\
 &  & \nopagebreak Reg. beds used  & $\text{11.970 $\pm$ 0.173}$ & $\text{12.463 $\pm$ 0.179}$ & $\text{13.056 $\pm$ 0.178}$ & $\text{13.875 $\pm$ 0.175}$ \\
 &  & \nopagebreak Beds reserved  & $\text{\phantom{0}1.154 $\pm$ 0.033}$ & $\text{\phantom{0}1.154 $\pm$ 0.033}$ & $\text{\phantom{0}1.163 $\pm$ 0.034}$ & $\text{\phantom{0}1.190 $\pm$ 0.036}$ \\
 &  & \nopagebreak Rooms added/removed  & $\text{\phantom{0}0.119 $\pm$ 0.004}$ & $\text{\phantom{0}0.119 $\pm$ 0.004}$ & $\text{\phantom{0}0.119 $\pm$ 0.004}$ & $\text{\phantom{0}0.120 $\pm$ 0.004}$ \\
 & \rule{0pt}{1.7\normalbaselineskip}ZGT & \nopagebreak Overbeds  & $\text{\phantom{0}0.366 $\pm$ 0.034}$ & $\text{\phantom{0}0.315 $\pm$ 0.034}$ & $\text{\phantom{0}0.269 $\pm$ 0.033}$ & $\text{\phantom{0}0.221 $\pm$ 0.033}$ \\
 &  & \nopagebreak Underbeds  & $\text{\phantom{0}3.952 $\pm$ 0.081}$ & $\text{\phantom{0}4.329 $\pm$ 0.088}$ & $\text{\phantom{0}4.779 $\pm$ 0.093}$ & $\text{\phantom{0}5.414 $\pm$ 0.101}$ \\
 &  & \nopagebreak Reg. beds used  & $\text{13.423 $\pm$ 0.135}$ & $\text{13.833 $\pm$ 0.132}$ & $\text{14.309 $\pm$ 0.128}$ & $\text{14.969 $\pm$ 0.119}$ \\
 &  & \nopagebreak Beds reserved  & $\text{\phantom{0}1.070 $\pm$ 0.030}$ & $\text{\phantom{0}1.036 $\pm$ 0.030}$ & $\text{\phantom{0}1.001 $\pm$ 0.029}$ & $\text{\phantom{0}0.976 $\pm$ 0.025}$ \\
 &  & \nopagebreak Rooms added/removed  & $\text{\phantom{0}0.092 $\pm$ 0.003}$ & $\text{\phantom{0}0.088 $\pm$ 0.003}$ & $\text{\phantom{0}0.085 $\pm$ 0.003}$ & $\text{\phantom{0}0.080 $\pm$ 0.003}$ \\
 & \rule{0pt}{1.7\normalbaselineskip}SKB & \nopagebreak Overbeds  & $\text{\phantom{0}0.255 $\pm$ 0.029}$ & $\text{\phantom{0}0.231 $\pm$ 0.029}$ & $\text{\phantom{0}0.207 $\pm$ 0.029}$ & $\text{\phantom{0}0.186 $\pm$ 0.028}$ \\
 &  & \nopagebreak Underbeds  & $\text{\phantom{0}2.805 $\pm$ 0.096}$ & $\text{\phantom{0}3.067 $\pm$ 0.101}$ & $\text{\phantom{0}3.388 $\pm$ 0.107}$ & $\text{\phantom{0}3.837 $\pm$ 0.119}$ \\
 &  & \nopagebreak Reg. beds used  & $\text{\phantom{0}6.561 $\pm$ 0.103}$ & $\text{\phantom{0}6.838 $\pm$ 0.101}$ & $\text{\phantom{0}7.174 $\pm$ 0.101}$ & $\text{\phantom{0}7.632 $\pm$ 0.103}$ \\
 &  & \nopagebreak Beds reserved  & $\text{\phantom{0}0.684 $\pm$ 0.025}$ & $\text{\phantom{0}0.685 $\pm$ 0.025}$ & $\text{\phantom{0}0.679 $\pm$ 0.023}$ & $\text{\phantom{0}0.675 $\pm$ 0.025}$ \\
 &  & \nopagebreak Rooms added/removed  & $\text{\phantom{0}0.061 $\pm$ 0.003}$ & $\text{\phantom{0}0.060 $\pm$ 0.003}$ & $\text{\phantom{0}0.060 $\pm$ 0.003}$ & $\text{\phantom{0}0.056 $\pm$ 0.003}$ \\
\rule{0pt}{1.7\normalbaselineskip}PU & \nopagebreak MST & \nopagebreak Overbeds  & $\text{\phantom{0}0.100 $\pm$ 0.011}$ & $\text{\phantom{0}0.099 $\pm$ 0.011}$ & $\text{\phantom{0}0.102 $\pm$ 0.012}$ & $\text{\phantom{0}0.096 $\pm$ 0.011}$ \\
 &  & \nopagebreak Underbeds  & $\text{\phantom{0}6.581 $\pm$ 0.073}$ & $\text{\phantom{0}6.582 $\pm$ 0.073}$ & $\text{\phantom{0}6.570 $\pm$ 0.073}$ & $\text{\phantom{0}6.614 $\pm$ 0.074}$ \\
 &  & \nopagebreak Reg. beds used  & $\text{22.916 $\pm$ 0.114}$ & $\text{22.917 $\pm$ 0.114}$ & $\text{22.903 $\pm$ 0.115}$ & $\text{22.943 $\pm$ 0.113}$ \\
 &  & \nopagebreak Beds reserved  & $\text{\phantom{0}0.825 $\pm$ 0.022}$ & $\text{\phantom{0}0.826 $\pm$ 0.022}$ & $\text{\phantom{0}0.835 $\pm$ 0.024}$ & $\text{\phantom{0}0.829 $\pm$ 0.023}$ \\
 &  & \nopagebreak Rooms added/removed  & $\text{\phantom{0}0.077 $\pm$ 0.003}$ & $\text{\phantom{0}0.077 $\pm$ 0.003}$ & $\text{\phantom{0}0.078 $\pm$ 0.003}$ & $\text{\phantom{0}0.077 $\pm$ 0.003}$ \\
 & \rule{0pt}{1.7\normalbaselineskip}ZGT & \nopagebreak Overbeds  & $\text{\phantom{0}0.187 $\pm$ 0.025}$ & $\text{\phantom{0}0.184 $\pm$ 0.025}$ & $\text{\phantom{0}0.181 $\pm$ 0.025}$ & $\text{\phantom{0}0.173 $\pm$ 0.024}$ \\
 &  & \nopagebreak Underbeds  & $\text{\phantom{0}3.313 $\pm$ 0.082}$ & $\text{\phantom{0}3.328 $\pm$ 0.082}$ & $\text{\phantom{0}3.359 $\pm$ 0.084}$ & $\text{\phantom{0}3.415 $\pm$ 0.086}$ \\
 &  & \nopagebreak Reg. beds used  & $\text{\phantom{0}9.072 $\pm$ 0.119}$ & $\text{\phantom{0}9.088 $\pm$ 0.119}$ & $\text{\phantom{0}9.126 $\pm$ 0.118}$ & $\text{\phantom{0}9.199 $\pm$ 0.118}$ \\
 &  & \nopagebreak Beds reserved  & $\text{\phantom{0}0.876 $\pm$ 0.026}$ & $\text{\phantom{0}0.882 $\pm$ 0.027}$ & $\text{\phantom{0}0.879 $\pm$ 0.026}$ & $\text{\phantom{0}0.895 $\pm$ 0.028}$ \\
 &  & \nopagebreak Rooms added/removed  & $\text{\phantom{0}0.089 $\pm$ 0.003}$ & $\text{\phantom{0}0.089 $\pm$ 0.003}$ & $\text{\phantom{0}0.089 $\pm$ 0.003}$ & $\text{\phantom{0}0.088 $\pm$ 0.003}$ \\
 & \rule{0pt}{1.7\normalbaselineskip}SKB & \nopagebreak Overbeds  & $\text{\phantom{0}0.277 $\pm$ 0.036}$ & $\text{\phantom{0}0.275 $\pm$ 0.036}$ & $\text{\phantom{0}0.268 $\pm$ 0.036}$ & $\text{\phantom{0}0.257 $\pm$ 0.035}$ \\
 &  & \nopagebreak Underbeds  & $\text{\phantom{0}1.908 $\pm$ 0.075}$ & $\text{\phantom{0}1.931 $\pm$ 0.076}$ & $\text{\phantom{0}1.986 $\pm$ 0.078}$ & $\text{\phantom{0}2.097 $\pm$ 0.080}$ \\
 &  & \nopagebreak Reg. beds used  & $\text{\phantom{0}4.754 $\pm$ 0.071}$ & $\text{\phantom{0}4.778 $\pm$ 0.070}$ & $\text{\phantom{0}4.836 $\pm$ 0.069}$ & $\text{\phantom{0}4.950 $\pm$ 0.066}$ \\
 &  & \nopagebreak Beds reserved  & $\text{\phantom{0}0.427 $\pm$ 0.017}$ & $\text{\phantom{0}0.432 $\pm$ 0.017}$ & $\text{\phantom{0}0.437 $\pm$ 0.019}$ & $\text{\phantom{0}0.434 $\pm$ 0.017}$ \\
 &  &\nopagebreak Rooms added/removed  & $\text{\phantom{0}0.040 $\pm$ 0.002}$ & $\text{\phantom{0}0.041 $\pm$ 0.002}$ & $\text{\phantom{0}0.040 $\pm$ 0.002}$ & $\text{\phantom{0}0.040 $\pm$ 0.002}$ \\
\hline 
\end{tabular}

\end{adjustbox}
\end{table}

\begin{table}
\small
\centering \caption{KPIs for SP-O for lookahead horizons 5, 4, and 3 days and when taking the $90\%$ CI upper bound for the predicted arrival intensity (UB arrival rate). }
\begin{adjustbox}{angle=90}
\centering
\small
\begin{tabular}{llcccc}
\hline
Hospital & KPI & Overbeds & UB arrival rate & Horizon 4 days & \multicolumn{1}{c}{Horizon 3 days} \\ 
\hline
\nopagebreak MST & \nopagebreak Overbeds  & $\text{\phantom{0}0.140 $\pm$ 0.018}$ & $\text{\phantom{0}0.105 $\pm$ 0.017}$ & $\text{\phantom{0}0.161 $\pm$ 0.020}$ & $\text{\phantom{0}0.240 $\pm$ 0.020}$ \\
 & \nopagebreak Underbeds  & $\text{\phantom{0}3.729 $\pm$ 0.104}$ & $\text{\phantom{0}4.403 $\pm$ 0.135}$ & $\text{\phantom{0}3.535 $\pm$ 0.106}$ & $\text{\phantom{0}3.132 $\pm$ 0.086}$ \\
 & \nopagebreak Reg. beds used  & $\text{12.187 $\pm$ 0.188}$ & $\text{13.189 $\pm$ 0.209}$ & $\text{11.822 $\pm$ 0.204}$ & $\text{11.260 $\pm$ 0.199}$ \\
 & \nopagebreak Beds reserved  & $\text{\phantom{0}0.828 $\pm$ 0.024}$ & $\text{\phantom{0}0.807 $\pm$ 0.024}$ & $\text{\phantom{0}0.916 $\pm$ 0.026}$ & $\text{\phantom{0}0.863 $\pm$ 0.023}$ \\
 & \nopagebreak Rooms added/removed  & $\text{\phantom{0}0.085 $\pm$ 0.003}$ & $\text{\phantom{0}0.073 $\pm$ 0.002}$ & $\text{\phantom{0}0.092 $\pm$ 0.003}$ & $\text{\phantom{0}0.088 $\pm$ 0.003}$ \\
\rule{0pt}{1.7\normalbaselineskip} ZGT & \nopagebreak Overbeds  & $\text{\phantom{0}0.096 $\pm$ 0.012}$ & $\text{\phantom{0}0.072 $\pm$ 0.012}$ & $\text{\phantom{0}0.112 $\pm$ 0.011}$ & $\text{\phantom{0}0.165 $\pm$ 0.015}$ \\
 & \nopagebreak Underbeds  & $\text{\phantom{0}3.397 $\pm$ 0.079}$ & $\text{\phantom{0}4.032 $\pm$ 0.086}$ & $\text{\phantom{0}2.952 $\pm$ 0.072}$ & $\text{\phantom{0}2.485 $\pm$ 0.068}$ \\
 & \nopagebreak Reg. beds used  & $\text{12.310 $\pm$ 0.171}$ & $\text{13.044 $\pm$ 0.190}$ & $\text{11.773 $\pm$ 0.179}$ & $\text{10.983 $\pm$ 0.187}$ \\
 & \nopagebreak Beds reserved  & $\text{\phantom{0}0.846 $\pm$ 0.024}$ & $\text{\phantom{0}0.774 $\pm$ 0.021}$ & $\text{\phantom{0}0.917 $\pm$ 0.026}$ & $\text{\phantom{0}0.865 $\pm$ 0.026}$ \\
 & \nopagebreak Rooms added/removed  & $\text{\phantom{0}0.084 $\pm$ 0.003}$ & $\text{\phantom{0}0.074 $\pm$ 0.002}$ & $\text{\phantom{0}0.091 $\pm$ 0.003}$ & $\text{\phantom{0}0.089 $\pm$ 0.003}$ \\
\rule{0pt}{1.7\normalbaselineskip} SKB & \nopagebreak Overbeds  & $\text{\phantom{0}0.090 $\pm$ 0.010}$ & $\text{\phantom{0}0.076 $\pm$ 0.011}$ & $\text{\phantom{0}0.107 $\pm$ 0.012}$ & $\text{\phantom{0}0.154 $\pm$ 0.012}$ \\
 & \nopagebreak Underbeds  & $\text{\phantom{0}2.509 $\pm$ 0.071}$ & $\text{\phantom{0}2.777 $\pm$ 0.084}$ & $\text{\phantom{0}2.345 $\pm$ 0.063}$ & $\text{\phantom{0}1.993 $\pm$ 0.052}$ \\
 & \nopagebreak Reg. beds used  & $\text{\phantom{0}7.221 $\pm$ 0.144}$ & $\text{\phantom{0}7.231 $\pm$ 0.149}$ & $\text{\phantom{0}7.220 $\pm$ 0.137}$ & $\text{\phantom{0}6.940 $\pm$ 0.139}$ \\
 & \nopagebreak Beds reserved  & $\text{\phantom{0}0.469 $\pm$ 0.018}$ & $\text{\phantom{0}0.442 $\pm$ 0.019}$ & $\text{\phantom{0}0.521 $\pm$ 0.022}$ & $\text{\phantom{0}0.470 $\pm$ 0.020}$ \\
 & \nopagebreak Rooms added/removed  & $\text{\phantom{0}0.042 $\pm$ 0.002}$ & $\text{\phantom{0}0.040 $\pm$ 0.002}$ & $\text{\phantom{0}0.045 $\pm$ 0.002}$ & $\text{\phantom{0}0.041 $\pm$ 0.002}$ \\
\rule{0pt}{1.7\normalbaselineskip} Region & \nopagebreak Overbeds  & $\text{\phantom{0}0.326 $\pm$ 0.035}$ & $\text{\phantom{0}0.253 $\pm$ 0.036}$ & $\text{\phantom{0}0.380 $\pm$ 0.038}$ & $\text{\phantom{0}0.559 $\pm$ 0.041}$ \\
 & \nopagebreak Underbeds  & $\text{\phantom{0}9.635 $\pm$ 0.181}$ & $\text{11.212 $\pm$ 0.210}$ & $\text{\phantom{0}8.832 $\pm$ 0.172}$ & $\text{\phantom{0}7.610 $\pm$ 0.148}$ \\
 & \nopagebreak Reg. beds used  & $\text{31.717 $\pm$ 0.240}$ & $\text{33.464 $\pm$ 0.240}$ & $\text{30.814 $\pm$ 0.239}$ & $\text{29.183 $\pm$ 0.242}$ \\
 & \nopagebreak Beds reserved  & $\text{\phantom{0}2.143 $\pm$ 0.042}$ & $\text{\phantom{0}2.023 $\pm$ 0.040}$ & $\text{\phantom{0}2.355 $\pm$ 0.043}$ & $\text{\phantom{0}2.198 $\pm$ 0.039}$ \\
 & \nopagebreak Rooms added/removed  & $\text{\phantom{0}0.211 $\pm$ 0.004}$ & $\text{\phantom{0}0.187 $\pm$ 0.004}$ & $\text{\phantom{0}0.228 $\pm$ 0.005}$ & $\text{\phantom{0}0.218 $\pm$ 0.004}$ \\
 & \rule{0pt}{1.7\normalbaselineskip}Avg. cost SP  & $\text{51.129 $\pm$ 1.407}$ & $\text{49.429 $\pm$ 1.459}$ & $\text{52.956 $\pm$ 1.533}$ & $\text{58.120 $\pm$ 1.676}$ \\
 & \nopagebreak Forecast cost SP (hor. 1)  & $\text{48.990 $\pm$ 1.238}$ & $\text{50.002 $\pm$ 1.423}$ & $\text{50.323 $\pm$ 1.293}$ & $\text{52.878 $\pm$ 1.328}$ \\
 & \nopagebreak Forecast cost SP (hor. 2)  & $\text{50.908 $\pm$ 1.325}$ & $\text{55.907 $\pm$ 1.644}$ & $\text{50.311 $\pm$ 1.356}$ & $\text{51.105 $\pm$ 1.307}$ \\
 & \nopagebreak Forecast cost SP (hor. 3)  & $\text{58.380 $\pm$ 1.797}$ & $\text{68.829 $\pm$ 2.325}$ & $\text{57.086 $\pm$ 1.867}$ & $\text{56.332 $\pm$ 1.738}$ \\
\hline 
\end{tabular}

\end{adjustbox}
\end{table}

\begin{table}
\small
\centering\caption{KPIs for SP-O, as well as the deterministic version of the stochastic programs, when taking the median and $85\%$ quantile over the arrivals and occupancy scenarios for each~day.}

\begin{adjustbox}{angle=90}
\centering
\small\begin{tabular}{llccc}
\hline
Hospital & KPI & Overbeds & Median Scen. & \multicolumn{1}{c}{Quantile Scen.} \\ 
\hline
\nopagebreak MST & \nopagebreak Overbeds  & $\text{\phantom{0}\phantom{0}0.140 $\pm$ 0.018}$ & $\text{\phantom{0}\phantom{0}0.831 $\pm$ 0.039}$ & $\text{\phantom{0}\phantom{0}0.107 $\pm$ 0.018}$ \\
 & \nopagebreak Underbeds  & $\text{\phantom{0}\phantom{0}3.729 $\pm$ 0.104}$ & $\text{\phantom{0}\phantom{0}1.718 $\pm$ 0.061}$ & $\text{\phantom{0}\phantom{0}5.512 $\pm$ 0.172}$ \\
 & \nopagebreak Reg. beds used  & $\text{\phantom{0}12.187 $\pm$ 0.188}$ & $\text{\phantom{0}\phantom{0}8.811 $\pm$ 0.179}$ & $\text{\phantom{0}14.439 $\pm$ 0.268}$ \\
 & \nopagebreak Beds reserved  & $\text{\phantom{0}\phantom{0}0.828 $\pm$ 0.024}$ & $\text{\phantom{0}\phantom{0}0.885 $\pm$ 0.019}$ & $\text{\phantom{0}\phantom{0}0.987 $\pm$ 0.034}$ \\
 & \nopagebreak Rooms added/removed  & $\text{\phantom{0}\phantom{0}0.085 $\pm$ 0.003}$ & $\text{\phantom{0}\phantom{0}0.104 $\pm$ 0.003}$ & $\text{\phantom{0}\phantom{0}0.086 $\pm$ 0.003}$ \\
\rule{0pt}{1.7\normalbaselineskip}ZGT & \nopagebreak Overbeds  & $\text{\phantom{0}\phantom{0}0.096 $\pm$ 0.012}$ & $\text{\phantom{0}\phantom{0}0.551 $\pm$ 0.028}$ & $\text{\phantom{0}\phantom{0}0.070 $\pm$ 0.012}$ \\
 & \nopagebreak Underbeds  & $\text{\phantom{0}\phantom{0}3.397 $\pm$ 0.079}$ & $\text{\phantom{0}\phantom{0}1.602 $\pm$ 0.054}$ & $\text{\phantom{0}\phantom{0}4.389 $\pm$ 0.120}$ \\
 & \nopagebreak Reg. beds used  & $\text{\phantom{0}12.310 $\pm$ 0.171}$ & $\text{\phantom{0}\phantom{0}9.000 $\pm$ 0.181}$ & $\text{\phantom{0}12.779 $\pm$ 0.235}$ \\
 & \nopagebreak Beds reserved  & $\text{\phantom{0}\phantom{0}0.846 $\pm$ 0.024}$ & $\text{\phantom{0}\phantom{0}0.880 $\pm$ 0.025}$ & $\text{\phantom{0}\phantom{0}0.892 $\pm$ 0.024}$ \\
 & \nopagebreak Rooms added/removed  & $\text{\phantom{0}\phantom{0}0.084 $\pm$ 0.003}$ & $\text{\phantom{0}\phantom{0}0.095 $\pm$ 0.003}$ & $\text{\phantom{0}\phantom{0}0.076 $\pm$ 0.003}$ \\
\rule{0pt}{1.7\normalbaselineskip}SKB & \nopagebreak Overbeds  & $\text{\phantom{0}\phantom{0}0.090 $\pm$ 0.010}$ & $\text{\phantom{0}\phantom{0}1.073 $\pm$ 0.047}$ & $\text{\phantom{0}\phantom{0}0.243 $\pm$ 0.034}$ \\
 & \nopagebreak Underbeds  & $\text{\phantom{0}\phantom{0}2.509 $\pm$ 0.071}$ & $\text{\phantom{0}\phantom{0}0.764 $\pm$ 0.039}$ & $\text{\phantom{0}\phantom{0}2.128 $\pm$ 0.080}$ \\
 & \nopagebreak Reg. beds used  & $\text{\phantom{0}\phantom{0}7.221 $\pm$ 0.144}$ & $\text{\phantom{0}\phantom{0}5.933 $\pm$ 0.148}$ & $\text{\phantom{0}\phantom{0}6.860 $\pm$ 0.168}$ \\
 & \nopagebreak Beds reserved  & $\text{\phantom{0}\phantom{0}0.469 $\pm$ 0.018}$ & $\text{\phantom{0}\phantom{0}0.444 $\pm$ 0.018}$ & $\text{\phantom{0}\phantom{0}0.460 $\pm$ 0.019}$ \\
 & \nopagebreak Rooms added/removed  & $\text{\phantom{0}\phantom{0}0.042 $\pm$ 0.002}$ & $\text{\phantom{0}\phantom{0}0.040 $\pm$ 0.002}$ & $\text{\phantom{0}\phantom{0}0.034 $\pm$ 0.002}$ \\
\rule{0pt}{1.7\normalbaselineskip}Region & \nopagebreak Overbeds  & $\text{\phantom{0}\phantom{0}0.326 $\pm$ 0.035}$ & $\text{\phantom{0}\phantom{0}2.455 $\pm$ 0.061}$ & $\text{\phantom{0}\phantom{0}0.419 $\pm$ 0.047}$ \\
 & \nopagebreak Underbeds  & $\text{\phantom{0}\phantom{0}9.635 $\pm$ 0.181}$ & $\text{\phantom{0}\phantom{0}4.084 $\pm$ 0.083}$ & $\text{\phantom{0}12.029 $\pm$ 0.199}$ \\
 & \nopagebreak Reg. beds used  & $\text{\phantom{0}31.717 $\pm$ 0.240}$ & $\text{\phantom{0}23.744 $\pm$ 0.245}$ & $\text{\phantom{0}34.079 $\pm$ 0.240}$ \\
 & \nopagebreak Beds reserved  & $\text{\phantom{0}\phantom{0}2.143 $\pm$ 0.042}$ & $\text{\phantom{0}\phantom{0}2.209 $\pm$ 0.034}$ & $\text{\phantom{0}\phantom{0}2.340 $\pm$ 0.044}$ \\
 & \nopagebreak Rooms added/removed  & $\text{\phantom{0}\phantom{0}0.211 $\pm$ 0.004}$ & $\text{\phantom{0}\phantom{0}0.239 $\pm$ 0.004}$ & $\text{\phantom{0}\phantom{0}0.196 $\pm$ 0.005}$ \\
 & \rule{0pt}{1.7\normalbaselineskip}Avg. cost SP  & $\text{\phantom{0}51.129 $\pm$ 1.407}$ & $\text{128.850 $\pm$ 2.542}$ & $\text{\phantom{0}57.299 $\pm$ 1.883}$ \\
 & \nopagebreak Forecast cost SP (hor. 1)  & $\text{\phantom{0}48.990 $\pm$ 1.238}$ & $\text{\phantom{0}85.416 $\pm$ 1.812}$ & $\text{\phantom{0}72.173 $\pm$ 2.988}$ \\
 & \nopagebreak Forecast cost SP (hor. 2)  & $\text{\phantom{0}50.908 $\pm$ 1.325}$ & $\text{\phantom{0}56.003 $\pm$ 1.257}$ & $\text{105.605 $\pm$ 4.429}$ \\
 & \nopagebreak Forecast cost SP (hor. 3)  & $\text{\phantom{0}58.380 $\pm$ 1.797}$ & $\text{\phantom{0}41.251 $\pm$ 1.445}$ & $\text{156.446 $\pm$ 5.864}$ \\
\hline 
\end{tabular}

\end{adjustbox}
\end{table}
\FloatBarrier

\end{document}